\documentclass[brochure,12pt]{bourbaki}
\pdfoutput=1
\usepackage{etex}
\usepackage{amssymb,amsfonts,amsmath,footnote}
\usepackage[utf8]{inputenc}
\usepackage[T1]{fontenc}
\usepackage[french]{babel}
\usepackage{lmodern}

\def\setminus{\mathchoice
    {\mathbin{\vrule height .92ex width 1.81ex depth -.58ex}}
    {\mathbin{\vrule height .92ex width 1.81ex depth -.58ex}}
    {\mathbin{\vrule height .65ex width 1.00ex depth -.43ex}}
    {\mathbin{\vrule height .50ex width 0.770ex depth -.34ex}}
}

\usepackage[
backend=biber,
style=authoryear, 
citestyle=authoryear-comp,
maxnames=4
]{biblatex}
\usepackage{csquotes}

\DeclareDelimFormat{nameyeardelim}{\addcomma\space}

\DeclareNameAlias{sortname}{given-family}

\renewcommand{\bibnamedash}{\leavevmode\raise3pt\hbox to3em{\hrulefill}\space}

\AtEveryBibitem{%
  \clearfield{issn} 
  \clearfield{doi} 

  \ifentrytype{online}{}{
    \clearfield{url}
  }
}

\renewbibmacro{in:}{%
    \ifentrytype{article}{}{\printtext{\bibstring{in}\intitlepunct}}}

\DeclareFieldFormat[article,periodical,inreference]{number}{\mkbibparens{#1}}
\DeclareFieldFormat[article,periodical,inreference]{volume}{\mkbibbold{#1}}
\renewbibmacro*{volume+number+eid}{%
    \printfield{volume}%
    \setunit*{\addthinspace}
    \printfield{number}%
    \setunit{\addcomma\space}%
    \printfield{eid}}

\DeclareFieldFormat[article,inbook,incollection]{title}{\enquote{#1}\addcomma} 

\date{Octobre 2017}
\bbkannee{70\`eme ann\'ee, 2017-2018}
\bbknumero{1138}
\title{GROUPES CONVEXES-COCOMPACTS EN RANG SUPÉRIEUR}
\subtitle{d'apr\`es Labourie, Kapovich, Leeb, Porti, \dots{}}
\author{Olivier GUICHARD}
\address{Universit\'e de Strasbourg, CNRS,
IRMA UMR 7501\\
F-67000, Strasbourg, France}
\email{olivier.guichard@math.unistra.fr}

\begin{document}
\maketitle

\noindent{\bf INTRODUCTION}

Les sous-groupes discrets des groupes de Lie, ou plutôt leurs actions sur les
espaces symétriques et plus généralement même les actions sur les immeubles
euclidiens, sont au cœur de cet exposé. Y seront présentés les sous-groupes
convexes-cocompacts d'isométries de l'espace hyperbolique $\mathbb{H}^n$ et
comment généraliser cette classe en rang supérieur. Spécifiquement, voici l'un
des résultats qui sera souligné et mis en contexte.

\begin{theo}[\cite{kapovichleebporti14_2}]
  \label{th1} Soit $\Gamma=\langle S\rangle$ un sous-groupe de
  $\mathbf{GL}_n(\mathbf{R})$ engendré par un ensemble fini $S$. On suppose
  \begin{description}
  \item[Régularité quasi-isométrique] Il existe $c>0$ et $C\geq 0$ tels que,
    pour tout $\gamma\in\Gamma$, $\mu_1(\gamma)-\mu_2(\gamma)\geq
    c\, |\gamma|_S-C$.
  \end{description}
  On a alors les conclusions suivantes
  \begin{description}
  \item[Hyperbolicité] le groupe $\Gamma$ est hyperbolique au sens de Gromov.
     Son bord de Gromov est noté $\partial_\infty \Gamma$.
  \item[Applications au bord] il existe des applications continues et
    $\Gamma$-équivariantes $\beta_1\colon \partial_\infty \Gamma \to
    \mathbb{P}^{n-1}(\mathbf{R})$ et $\beta_{n-1}\colon \partial_\infty \Gamma \to
    \mathbb{P}^{n-1*}(\mathbf{R})= \mathrm{Gr}_{n-1}(\mathbf{R}^n)$ telles que
    \[
      \begin{array}{ll}
        \forall t\in \partial_\infty \Gamma, & \beta_1(t)\subset
                                               \beta_{n-1}(t)\\
        \forall t\neq t'\in \partial_\infty \Gamma, & \beta_1(t) \cap
                                               \beta_{n-1}(t')=0.
      \end{array}
    \] 
  \item[Contraction exponentielle] Il existe $k>0$ et $K\geq 0$ tels que, pour
    tout rayon géodésique $(\gamma_p)_{p\in\mathbf{N}}$ dans le graphe de Cayley de
    $\Gamma$ avec $\gamma_0= e_\Gamma$ dont le point limite dans
    $\partial_\infty \Gamma$ est noté $\gamma_\infty$, on a, pour tout $p\in
    \mathbf{N}$,
    \[ \log \bigl\| \gamma_{p}^{-1}|_{T_{\beta_1(\gamma_\infty)}
        \mathbb{P}^{n-1}(\mathbf{R})}\bigr\| \geq k\, p -K=k\, |\gamma_p|_S-K. \]
  \end{description}
Réciproquement, un sous-groupe $\Gamma$ vérifiant ces conclusions satisfait à
l'hypothèse de régularité quasi-isométrique.
\end{theo}

Précisons d'abord quelques notations utilisées dans cet énoncé.

Lorsqu'un groupe $\Gamma$ a une partie génératrice $S$ donnée, la \emph{fonction
longueur} associée est notée $|\cdot|_S\colon \Gamma \to \mathbf{N}$, explicitement,
pour tout $\gamma\in\Gamma$, $|\gamma|_S = \inf \{ p\in \mathbf{N} \mid
\exists \gamma_1, \dots, \gamma_p \in S\cup S^{-1},\, \gamma = \gamma_1 \cdots
\gamma_p\}$ (le produit vide étant égal à $e_\Gamma$). Une suite
$(\gamma_p)_{p\in \mathbf{N}}$ de $\Gamma$ est dès lors une géodésique dans le
graphe de Cayley si, pour tout $p\in\mathbf{N}$, $\gamma_{p}^{-1} \gamma_{p+1}
\in S\cup S^{-1}$ et $|\gamma_{0}^{-1}\gamma_p|_S=p$.

Le bord de Gromov de $\Gamma$ est l'ensemble des classes d'équivalence de rayons
géodésiques pour la relation \og être à distance de Hausdorff bornée\fg{}.

Si $S_1$ est une seconde partie génératrice \emph{finie} de $\Gamma$, il
existe $c_1>0$ tel que, pour tout $\gamma\in\Gamma$,
$|\gamma|_{S_1} \geq c_1 \, |\gamma|_S$. Ainsi, le choix de la
partie génératrice finie $S$ dans l'énoncé ci-dessus importe donc peu et la
notation $|\cdot|_\Gamma$ sera dorénavant adoptée pour désigner l'une des
fonctions longueurs associées à une partie génératrice finie.

Si $g$ est un élément de $\mathbf{GL}_n(\mathbf{R})$, $\mu_1(g), \dots,
\mu_n(g)$ désignent les logarithmes des valeurs principales de $g$,
c'est-à-dire $e^{2\mu_1(g)}, \dots, e^{2\mu_n(g)}$ est la liste décroissante
des valeurs propres de $g \, {}^t\! g$. Dans la conclusion de l'énoncé, $g|_{T_x
\mathbb{P}^{n-1}(\mathbf{R})}$ désigne l'application tangente au point $x$ du
difféomorphisme de $\mathbb{P}^{n-1}(\mathbf{R})$ induit par l'élément $g$ de
$\mathbf{GL}_n(\mathbf{R})$, c'est une application linéaire de $T_x
\mathbb{P}^{n-1}(\mathbf{R})$ dans $T_{g\cdot x}
\mathbb{P}^{n-1}(\mathbf{R})$. Sa norme $\| g|_{T_x
\mathbb{P}^{n-1}(\mathbf{R})}\|$ est la norme d'opérateur calculée pour une
métrique riemannienne sur $\mathbb{P}^{n-1}(\mathbf{R})$. La compacité de
$\mathbb{P}^{n-1}(\mathbf{R})$ implique que le choix de cette structure
riemannienne n'influe que sur les constantes $k$ et $K$ de la conclusion du
théorème.

\medskip

Le théorème énoncé ici est un \emph{cas particulier} (emblématique, bien sûr)
des résultats de Kapovich, Leeb et Porti. En effet ceux-ci autorisent :
\begin{itemize}
\item plus de généralités sur la \emph{source}. Il n'est pas nécessaire de
  travailler avec un sous-groupe de type fini $\Gamma$ mais une application
  $\varphi\colon Z\to \mathbf{GL}_n(\mathbf{R})$ d'un espace métrique $(Z,d_Z)$
  géodésique localement compact dans $\mathbf{GL}_n(\mathbf{R})$, l'hypothèse
  s'écrivant ici : $\forall z_1,z_2 \in Z$, $\mu_1\bigl( \varphi(z_1)^{-1}
  \varphi(z_2)\bigr) - \mu_2\bigl( \varphi(z_1)^{-1} \varphi(z_2)\bigr) \geq c\,
  d_Z(z_1,z_2)-C$. La conclusion affirme que $Z$ est hyperbolique au sens de
  Gromov et qu'il existe des applications continues $\beta_1\colon \partial_\infty
  Z \to \mathbb{P}^{n-1}(\mathbf{R})$, $\beta_{n-1}\colon \partial_\infty
  Z \to \mathbb{P}^{n-1*}(\mathbf{R})$ satisfaisant des propriétés de
  transversalité et de contraction similaires à celles énoncées plus haut.
\item plus de généralités sur le \emph{but}. Le même résultat est
  valable en remplaçant $\mathbf{R}$ par $\mathbf{C}$ ou par une extension
  finie de $\mathbf{Q}_p$. En fait, on peut aussi considérer des sous-groupes
  discrets d'un groupe algébrique $\mathbf{G}$ semi-simple sur un corps topologique
  localement compact et même sur un corps non-archimédien à valuation non discrète. L'hypothèse est alors que les projections de Cartan (qui
  sont des éléments de la chambre de Weyl fermée) des éléments de $\Gamma$
  sont quantitativement \og loin\fg{} d'une réunion de murs de la chambre de
  Weyl\footnote{Les énoncés des articles de Kapovich, Leeb et Porti adoptent
    un point de vue \og dual\fg{} en demandant que les projections de Cartan
    soient quantitativement \og proches\fg{} de la face de la chambre de Weyl
    égale à l'intersection des murs restants.}. La conclusion est encore
  l'hyperbolicité au sens de Gromov de $\Gamma$ et l'existence d'une
  application, du bord de Gromov~$\partial_\infty \Gamma$ dans la variété
  drapeau associée à cette réunion finie de murs, satisfaisant ces propriétés
  de transversalité et de contraction.

  {
    Pour le cas de $\mathbf{GL}_n(\mathbf{R})$, la chambre de Weyl fermée $A^+$
  est l'ensemble des matrices diagonales à coefficients décroissants, la
  projection d'un élément $g$ de $\mathbf{GL}_n(\mathbf{R})$ est la matrice
  dont la diagonale est $(\mu_1(g), \mu_2(g),\dots, \mu_n(g))$. Les murs
  correspondant à l'énoncé du théorème sont $\{ (m_1, m_2, \dots,
  m_n)\in A^+ \mid m_1=m_2\}$ et $\{ (m_1, \dots, m_n)\in A^+ \mid
  m_{n-1}=m_n\}$, 
  l'éloignement de ce second mur provient de l'hypothèse du théorème appliquée
  à $\gamma^{-1}$. La variété drapeau associée à ces deux murs est $\mathcal 
    {F}_{1,n-1} = \{ (\ell,h) \in \mathbb{P}^{n-1}(\mathbf{R}) \times
  \mathbb{P}^{n-1*}(\mathbf{R}) \mid \ell\subset h\}$, les
  applications~$\beta_1$ 
  et~$\beta_{n-1}$ se combinant en une application continue, équivariante et
  transverse $\beta\colon \partial_\infty \Gamma \to {\mathcal
    F}_{1,n-1}$. Notons que, même pour $\mathbf{GL}_n( \mathbf{R})$, le
  théorème principal est aussi présenté dans une généralité restreinte
  puisque tous les autres choix de murs sont possibles conduisant à des
  applications au bord dans des variétés drapeaux différentes. \par }

  La condition sur les projections de Cartan peut également s'exprimer en
  termes de l'action sur l'espace symétrique, ou bien sur l'immeuble de Bruhat
  et Tits dans le cas où le groupe algébrique est défini sur un corps
  ultramétrique, associé à $\mathbf{G}$. La position relative de deux points
  $x$ et $y$   de cet espace symétrique ou immeuble euclidien est encore un
  élément de la chambre de Weyl du groupe $\mathbf{G}$ et la projection de
  Cartan d'un élément $g\in\mathbf{G}$ est égale à la position relative de
  $x_0$ et $g\cdot x_0$ (où $x_0$ est un point base de l'espace symétrique ou
  immeuble euclidien en question). Ce point de vue permet donc de traiter
  d'une même manière ces deux cas mais également d'inclure le cas de
  sous-groupes du groupe des isométries d'un immeuble euclidien non
  nécessairement localement compact. La démonstration du théorème fait
  apparaître ce type d'immeubles comme cônes asymptotiques et il est naturel
  de  les inclure dans la discussion.

  Pour limiter le nombre de notions à introduire, cet exposé n'adoptera que
  peu ou prou ce langage géométrique au prix de présenter des résultats
  peut-être partiels. Ce point de vue est largement développé dans les
  articles de 
  \textcite{kapovichleeb15,kapovichleebporti14_2,kapovichleebportiv2,kapovichleebporti17} ainsi que dans les notes de
  cours et panorama qu'ils ont écrits \parencite{kapovichleeb17,kapovichleebporti15}. Notons cependant que même
  pour $n=3$ et pour $\Gamma \simeq \mathbb{F}_2$ le groupe libre à deux
  générateurs (auquel cas l'hyperbolicité de~$\Gamma$ n'est plus à établir) le
  résultat est nouveau. 
\end{itemize}

\smallskip

Le but de ce texte est de donner des éléments de la démonstration du théorème~\ref{th1},
de donner plusieurs autres caractérisations de la classe des sous-groupes
(nécessairement discrets) qui apparaissent ici, de faire le lien avec la
classe des sous-groupes convexes-cocompacts, d'énumérer certaines propriétés
géométriques et dynamiques remarquables de cette classe de sous-groupes.

\smallskip

Aussi bien les hypothèses que les conclusions de ce théorème impliquent
uniquement la projection de $\Gamma$ dans $\mathbf{PGL}_n(\mathbf{R})$. Comme
il est parfois commode de se débarrasser du centre de
$\mathbf{GL}_n(\mathbf{R})$, la suite de cet exposé privilégiera de temps en
temps les sous-groupes de $\mathbf{PGL}_n(\mathbf{R})$, voire de
$\mathbf{SL}_n(\mathbf{R})$. Nous ne l'avons pas écrit mais l'action de
$\Gamma$ sur $\mathbb{P}^{n-1*}(\mathbf{R})$ vérifie des propriétés de
contraction équivalentes.

Il peut être aussi utile d'observer que les propriétés énoncées dans ce
théorème sont \og stables par passage à un sous-groupe d'indice fini\fg{};
précisément
\begin{itemize}
\item si l'hypothèse de régularité quasi-isométrique est satisfaite pour un
  sous-groupe d'indice fini de $\Gamma$, alors elle l'est pour $\Gamma$
  (quitte à changer les constantes $c$ et $C$);
\item de même, si les conclusions sont satisfaites  pour un
  sous-groupe d'indice fini de~$\Gamma$, alors elles le sont pour $\Gamma$
  (quitte à changer les constantes $k$ et $K$).
\end{itemize}

\section{Groupes convexes-cocompacts d'isométries de l'espace hyperbolique}
\label{sec:group-conv-cocomp}

Toujours pour rester concret, la discussion suivante est restreinte aux
sous-groupes d'isométries de l'espace hyperbolique réel $\mathbb{H}^n$, $n\geq 2$,
mais son contexte naturel serait plutôt celui des groupes d'isométries de
variétés de Cartan et Hadamard de courbure sectionnelle majorée par une
constante strictement négative.

\subsection{L'espace hyperbolique et ses convexes}
\label{sec:lesp-hyperb-et}

Le sous-groupe des matrices orthogonales pour $Q_{n,1} =
\bigl(\begin{smallmatrix}
  1_n & 0 \\ 0 & -1
\end{smallmatrix}\bigr)
$ est noté $\mathbf{O}(n,1) = \{ g\in \mathbf{GL}_{n+1}(\mathbf{R}) \mid {}^t\!
g Q_{n,1} g = Q_{n,1}\}$. La forme quadratique associée à $Q_{n,1}$ est
$q_{n,1}\colon \mathbf{R}^{n+1}\to \mathbf{R} / x \mapsto {}^t\! x Q_{n,1}x$. 
Le modèle projectif de l'espace hyperbolique est
\[ \mathbb{H}^n = \{ \ell\in \mathbb{P}^n(\mathbf{R}) \mid q_{n,1}|_{ \ell\setminus
    \{0\}} <0\}\ ;\]
il s'identifie à l'hyperboloïde $\{ x= (x_1, \dots , x_{n+1})
\in \mathbf{R}^{n+1} \mid q_{n,1}(x)=-1,\ x_{n+1}>0\}$. Son adhérence dans
$\mathbb{P}^n(\mathbf{R})$ est $\overline{\mathbb{H}}{}^n =\{ \ell\in
\mathbb{P}^n(\mathbf{R}) \mid q_{n,1}|_\ell \leq 0\}$ et est la réunion de
$\mathbb{H}^n$ et de $\partial_\infty \mathbb{H}^n=\{ \ell\in
\mathbb{P}^n(\mathbf{R}) \mid q_{n,1}|_\ell = 0\}$. Aussi bien $\mathbb{H}^n$ que
son adhérence sont contenus dans la carte affine $\mathbb{A}^n =
\mathbb{P}^n(\mathbf{R}) \setminus \mathbb{P}(\{ x_{n+1}=0\})$. Dans les
coordonnées naturelles sur cette carte affine, $\mathbb{H}^n$ et
$\overline{\mathbb{H}}{}^n$ sont respectivement la boule unité ouverte et la
boule unité fermée.

L'action de $\mathbf{O}(n,1)$ sur $\mathbb{H}^n$ est transitive et le
stabilisateur du point $\ell_0 = \mathbf{R} e_{n+1}$ ($(e_i)_{i=1,\dots,n+1}$
désigne la base canonique de $\mathbf{R}^{n+1}$) est $\mathbf{O}(n)\times
\mathbf{O}(1)$. Puisque ce stabilisateur est compact, il existe une métrique
riemannienne $\mathbf{O}(n,1)$-invariante sur~$\mathbb{H}^n$. La distance
associée $d_{\mathbb{H}^n}$ s'exprime aisément en termes de la géométrie
projective : si~$\ell$ et~$\ell'$ sont dans~$\mathbb{H}^n$, alors $d_{\mathbb{H}^n}(\ell,\ell') =
1/2 \, \bigl| \log | [\ell,\ell';b,b']|\bigr|$ où $b$, $b'$ sont deux éléments
distincts de~$\partial_\infty \mathbb{H}^n$ tels que $\ell,\ell',b,b'$ appartiennent
à une même droite projective (si $\ell\neq \ell'$, alors $\{ b,b'\}$ est
l'intersection de $\partial_\infty \mathbb{H}^n$ et de la droite projective
engendrée par $\ell$ et $\ell'$) et où $[\ell,\ell';b,b'] = \frac{ \ell-b}{\ell-b'} \frac{
  \ell'-b'}{\ell'-b}$ est le birapport des quatre points $\ell,\ell',b,b'$ calculé dans
une identification quelconque de la droite projective contenant ces points
avec $\mathbb{P}^1(\mathbf{R})$.

Le groupe des isométries de $(\mathbb{H}^n, d_{\mathbb{H}^n})$ est alors
$\mathbf{PO}(n,1) = \mathbf{O}(n,1)/ \{ \pm 1_{n+1}\}$.

Deux points de $\mathbb{H}^n$ sont reliés par un unique segment géodésique qui
est égal au segment affine les reliant dans $\mathbb{A}^n$. Plus généralement,
un point de $\mathbb{H}^n$ et un point de $\partial_\infty \mathbb{H}^n$
(respectivement deux points distincts  de $\partial_\infty \mathbb{H}^n$) sont
reliés par un unique rayon géodésique (respectivement une unique géodésique)
égal à un segment affine semi-ouvert (respectivement égale à un segment affine
ouvert).

Un sous-ensemble $C$ de $\mathbb{H}^n$ ou de $\overline{\mathbb{H}}{}^n$ est
dit \emph{convexe} s'il est convexe au sens de la géométrie affine dans la
carte affine $\mathbb{A}^n$. Une condition équivalente, si $C$ est inclus dans
$\mathbb{H}^n$, est de demander que $C$ est géodésiquement convexe.

\subsection{Groupes convexes-cocompacts, leurs ensembles limites}
\label{sec:group-conv-cocomp-1}

\begin{defi}
  \label{defi:group-conv-cocomp}
  Un sous-groupe $\Gamma$ de $\mathbf{O}(n,1)$ est dit \emph{convexe-cocompact
  } s'il existe un sous-ensemble non vide $C$ de $\mathbb{H}^n$, convexe,
  $\Gamma$-invariant et sur lequel l'action de $\Gamma$ est propre et cocompacte.
\end{defi}
Le convexe $C$ est alors nécessairement fermé et le sous-groupe $\Gamma$ est
nécessairement discret. En remplaçant éventuellement $C$ par un voisinage
métrique $N_\delta(C) = \{ \ell\in \mathbb{H}^n \mid
\mathrm{dist}_{\mathbb{H}^n}(\ell,C)\leq \delta\}$, le quotient $\Gamma
\backslash C$ est une orbivariété à bord dont le groupe fondamental orbifold
est égal à $\Gamma$. Le groupe $\Gamma$ est donc de \emph{type fini} et même
de présentation finie. Les sous-groupes finis de $\mathbf{O}(n,1)$ sont
toujours convexes-cocompacts ; il est alors sensé de n'étudier que des groupes
infinis. Un sous-groupe cyclique infini $\langle \gamma\rangle$ de
$\mathbf{O}(n,1)$ est convexe-cocompact si et seulement si l'élément $\gamma$
est hyperbolique ; il a alors un unique axe de translation qui est une
géodésique de $\mathbb{H}^n$, le convexe $C$ peut être choisi égal à cet axe.

Si $\Gamma$ est un sous-groupe de $\mathbf{O}(n,1)$, son \emph{ensemble
  limite} $\Lambda_\Gamma \subset \partial_\infty \mathbb{H}^n$ est défini
comme l'ensemble des points d'accumulation de l'orbite $\Gamma\cdot \ell_0
\subset \mathbb{H}^n$ ($\ell_0=\mathbf{R} e_{n+1} $) dans le bord à l'infini $\partial_\infty
\mathbb{H}^n$. C'est aussi l'ensemble de points d'accumulation de toute orbite
de $\Gamma$ dans~$\mathbb{H}^n$. L'ensemble limite est vide si et seulement si
le sous-groupe $\Gamma$ est borné. L'ensemble limite est réduit à un seul élément
$x\in \partial_\infty \mathbb{H}^n$ si et seulement si $\Gamma$ est
infini et est contenu dans le sous-groupe $\{ g\in \mathbf{O}(n,1) \mid g|_x =
\pm \mathrm{Id}|_x\} \subset \mathrm{Stab}_{\mathbf{O}(n,1)}(x)$. Lorsque~$\Gamma$ laisse invariant un sous-ensemble convexe $C$ de $\mathbb{H}^n$,
l'ensemble limite $\Lambda_\Gamma$ est alors inclus dans le bord idéal $\partial_\infty C =
\overline{C} \cap \partial_\infty \mathbb{H}^n$ ($\overline{C}$ est
l'adhérence de $C$ dans $\mathbb{P}^n(\mathbf{R})$). Si, de plus, le groupe~$\Gamma$ est
convexe-cocompact, l'inclusion inverse est vérifiée et il y a l'égalité
$\Lambda_\Gamma = \partial_\infty C$.

Ainsi, quand $\Gamma$ est convexe-cocompact et infini, il est alors
impossible d'avoir $\sharp \Lambda_\Gamma =1$ (le feuilletage horocyclique
centré en $\Lambda_\Gamma$ serait $\Gamma$-invariant), de sorte que l'on peut toujours
considérer l'enveloppe convexe $\mathrm{Conv}(\Lambda_\Gamma)$ de
$\Lambda_\Gamma$ dans $\mathbb{A}^n$ et que
$\mathrm{Conv}_{\mathbb{H}^n}(\Lambda_\Gamma)$, son intersection avec
$\mathbb{H}^n$, est non vide. Ce convexe
$\mathrm{Conv}_{\mathbb{H}^n}(\Lambda_\Gamma)$ est alors fermé,
\mbox{$\Gamma$-invariant} et inclus dans $C$ ; de là l'action de $\Gamma$ sur
$\mathrm{Conv}_{\mathbb{H}^n}(\Lambda_\Gamma)$ est propre et cocompacte. Le
cas où $\Lambda_\Gamma = \{ x,x'\}$ ($x\neq x'$) est aussi particulier : dans
ce cas, un sous-groupe d'indice fini de $\Gamma$ est infini cyclique engendré
par un élément hyperbolique. On appelle un groupe~$\Gamma$ \emph{élémentaire} si
$\sharp \Lambda_\Gamma =0,1$ ou 2 ; sinon $\Lambda_\Gamma$ a la puissance du
continu et peut être, par exemple, homéomorphe à l'ensemble de Cantor.

\subsection{Caractérisation métrique}
\label{sec:caract-metr}

Les valeurs principales d'un élément $g$ de $\mathbf{O}(n,1)$ sont $(
\mu_1(g), \dots, \mu_{n+1}(g)) = ( \mu(g), 0, \dots, 0, -\mu(g))$ où
$\mu(g)\in \mathbf{R}_+$. La double classe $\bigl( \mathbf{O}(n) \times
\mathbf{O}(1)\bigr) g \bigl( \mathbf{O}(n) \times \mathbf{O}(1)\bigr)$
contient un unique élément de la forme $\Bigl\{
\Bigl(\begin{smallmatrix}
  1_{n-1} & & \\ & \cosh t & \sinh t \\ & \sinh t & \cosh t
\end{smallmatrix}\Bigr),\, t\in \mathbf{R}_+ \Bigr\}
$ pour $t=\mu(g)$. La distance $d_{\mathbb{H}^n}( \ell_0, g\cdot \ell_0)$, toujours
avec $\ell_0= \mathbf{R}e_{n+1}$, est égale à $\mu(g)$.

On a alors la caractérisation métrique suivante.

\begin{prop}\label{prop:conv-cocomp-quasi-iso}
  Soit $\Gamma$ un sous-groupe de $\mathbf{O}(n,1)$. Le groupe $\Gamma$ est
  convexe-cocompact si et seulement si $\Gamma$ est de type fini et il existe
  $c>0$, $C\geq 0$ tels que, pour tout $\gamma\in\Gamma$, $\mu(\gamma) \geq c
  \, |\gamma|_\Gamma  -C$.
\end{prop}

On formule souvent cette conclusion en disant que le groupe $\Gamma$ est
\emph{quasi-isométriquement plongé} dans $\mathbf{O}(n,1)$. Le sens direct de cette
proposition suit de ce que l'on a dit plus haut et du lemme de Milnor et \v{S}varc
tandis que le sens réciproque utilise deux ingrédients : le fait que les
quasi-géodésiques de $\mathbb{H}^n$ sont universellement proches de
géodésiques (\og lemme de Morse\fg{}) et le fait que les simplexes géodésiques
sont universellement fins.

\subsection{Caractérisation dynamique}
\label{sec:caract-dynam}

Le flot géodésique $(\varphi_t)_{t\in\mathbf{R}}$ sur le fibré unitaire
tangent $T^1 \mathbb{H}^n$ commute avec l'action de $\mathbf{O}(n,1)$. Si
$\Gamma$ est un sous-groupe de $\mathbf{O}(n,1)$, on définit l'\emph{ensemble
récurrent} $\mathcal{E}_\Gamma \subset T^1 \mathbb{H}^n$ de la manière suivante
: un vecteur $v$ de $T^1 \mathbb{H}^n$ appartient à $\mathcal{E}_\Gamma$ si et
seulement si, pour tout voisinage $U$ de $v$ dans $T^1 \mathbb{H}^n$, il
existe une suite $(\gamma_k,t_k)_{k\in \mathbf{Z}}$ de $\Gamma \times
\mathbf{R}$ avec la propriété que, pour tout $k\in \mathbf{Z}$, $\gamma_k \cdot
\varphi_{t_k}(v) = \varphi_{t_k}(\gamma_k \cdot v )$ appartient à $U$ et
$\lim_{k\to \pm \infty}
t_k = \pm\infty$. 
L'ensemble récurrent $\mathcal{E}_\Gamma$ est fermé,
$\Gamma$-invariant et $\varphi_t$-invariant. Si $\Gamma$ est discret, le
quotient $\Gamma \backslash \mathcal{E}_\Gamma$ est l'ensemble récurrent pour
le flot géodésique sur le fibré unitaire tangent $\Gamma\backslash T^1
\mathbb{H}^n$ de l'orbivariété riemannienne $\Gamma\backslash  \mathbb{H}^n$.

Voici une caractérisation dynamique des sous-groupes convexes-cocompacts.

\begin{prop}\label{prop:conv-cocomp-recurrent}
  Soit $\Gamma$ un sous-groupe de $\mathbf{O}(n,1)$. Le groupe $\Gamma$ est
  convexe-cocompact si et seulement si l'action de $\Gamma$ sur
  $\mathcal{E}_\Gamma$ est propre et cocompacte.
\end{prop}
{
\noindent Lorsque $\Gamma$ est convexe-cocompact, ou plus généralement s'il existe un
convexe fermé $C\subset \mathbb{H}^n$ invariant par $\Gamma$, alors l'ensemble
$\mathcal{E}_\Gamma$ est inclus dans l'ensemble des vecteurs $v\in T^1
\mathbb{H}^n$ tels que, pour tout $t\in \mathbf{R}$, le point base de
$\varphi_t(v)$ appartient à $C$, i.e.\ $\mathcal{E}_\Gamma \subset T^1
\mathbb{H}^n |_C$. Quand $\Gamma$ est convexe-cocompact, $\Gamma\backslash T^1
\mathbb{H}^n |_C$ est compact et $\Gamma \backslash \mathcal{E}_\Gamma$
aussi. Réciproquement, si $\Gamma \backslash \mathcal{E}_\Gamma$ est compact,
le groupe $\Gamma$ est alors de type fini et l'injection de $\Gamma$ dans
$\mathbf{O}(n,1)$ est un plongement quasi-isométrique. La
proposition~\ref{prop:conv-cocomp-quasi-iso} s'applique.\par
}

\subsection{Caractérisations par l'action à l'infini de l'espace hyperbolique}
\label{sec:caract-par-lact}

La convexe-cocompacité peut aussi être caractérisée par l'action de $\Gamma$
sur $\overline{ \mathbb{H}}{}^n = \mathbb{H}^n \cup \partial_\infty
\mathbb{H}^n$.
\begin{prop}
  \label{prop:conv-cocomp-hnbar}
  Soit $\Gamma$ un sous-groupe de $\mathbf{O}(n,1)$. Le groupe $\Gamma$ est
  convexe-cocompact si et seulement si l'action de $\Gamma$ sur
  $\overline{ \mathbb{H}}{}^n \setminus \Lambda_\Gamma$ est propre et cocompacte.
\end{prop}
{
  \noindent En effet si $\Gamma$ est convexe-cocompact, la projection $\mathbb{H}^n\to C$
s'étend à $\overline{ \mathbb{H}}{}^n \setminus \partial_\infty C = \overline{
  \mathbb{H}}{}^n \setminus \Lambda_\Gamma$ et cette extension est propre et
$\Gamma$-équivariante ; de là suivent la propreté et la cocompacité de
l'action de $\Gamma$ sur $\overline{ \mathbb{H}}{}^n \setminus
\Lambda_\Gamma$. Réciproquement supposons l'action de~$\Gamma$ sur $\overline{
  \mathbb{H}}{}^n \setminus \Lambda_\Gamma$ propre et cocompacte. Si
$\Lambda_\Gamma=\emptyset$, alors $\Gamma$ est fini et le cas $\sharp
\Lambda_\Gamma =1$ est impossible (existence du feuilletage horocyclique). Dans les autres cas, le
convexe $\mathrm{Conv}_{\mathbb{H}^n} ( \Lambda_\Gamma)$ est inclus et fermé
dans $\overline{ \mathbb{H}}{}^n \setminus \Lambda_\Gamma$, le groupe $\Gamma$
agit donc proprement sur  $\mathrm{Conv}_{\mathbb{H}^n} ( \Lambda_\Gamma)$ et
$\Gamma \backslash  \mathrm{Conv}_{\mathbb{H}^n} ( \Lambda_\Gamma)$ est fermé
dans $\Gamma \backslash ( \overline{ \mathbb{H}}{}^n \setminus
\Lambda_\Gamma)$ donc est compact.\par
}

En particulier l'action de $\Gamma$ sur $\partial_\infty \mathbb{H}^n
\setminus \Lambda_\Gamma$ est propre et cocompacte. Il faut se garder de croire
que la réciproque est exacte, par exemple les travaux de 
\textcite{bers-bdry-teich} impliquent l'existence de sous-groupes $\Gamma$ de
$\mathbf{O}(3,1)$, isomorphes au groupe fondamental d'une surface compacte, pour
lesquels $\Gamma \backslash ( \partial_\infty \mathbb{H}^3 \setminus
\Lambda_\Gamma)$ est compact mais $\Gamma \backslash (  \overline{\mathbb{H}}{}^3 \setminus
\Lambda_\Gamma)$ ne l'est pas. Néanmoins la caractérisation suivante
n'implique que l'action de $\Gamma$ sur $\partial_\infty \mathbb{H}^n$.

\medskip

Soit $d_{\partial_\infty \mathbb{H}^n}$ une distance sur $\partial_\infty
\mathbb{H}^n$ provenant d'une métrique riemannienne ($\partial_\infty
\mathbb{H}^n$ étant une variété compacte, le choix particulier de la métrique
riemannienne n'influe que sur certaines constantes dans la suite). L'action de
$\Gamma$ est dite \emph{dilatante} au point~$z$ de~$\partial_\infty
\mathbb{H}^n$ s'il existe $\gamma\in\Gamma$, $c>1$ et $U\subset \partial_\infty
\mathbb{H}^n$ un voisinage de $z$ tels que, pour tous~$z_1$, $z_2$ dans $U$,
$d_{\partial_\infty \mathbb{H}^n} ( \gamma\cdot z_1, \gamma\cdot z_2) \geq c\, 
d_{\partial_\infty \mathbb{H}^n}(z_1,z_2)$.

\begin{theo}
  \label{theo:group-conv-cocomp-dilatante-dHn}
  Soit $\Gamma$ un sous-groupe discret de $\mathbf{O}(n,1)$. Le groupe
  $\Gamma$ est convexe-cocompact si et seulement si, pour tout $z\in
  \Lambda_\Gamma$, l'action de $\Gamma$ est dilatante (dans $\partial_\infty
\mathbb{H}^n$) au point $z$.
\end{theo}
La démonstration de ce résultat est un peu plus délicate que les
caractérisations précédentes ; c'est pourquoi il se trouve ici sous
l'intitulé de théorème. Il faut en fait faire le lien entre l'action au bord
et l'action sur $\mathbb{H}^n$. Il est également possible de détecter, pour un
sous-groupe $\Gamma$ donné, si un point $z$ de $\Lambda_\Gamma$ est dilatant
selon la façon dont l'orbite $\Gamma\cdot \ell_0$ s'accumule sur $z$ (de manière
\og transverse\fg{}/\og conique\fg{} ou non par rapport au bord
$\partial_\infty \mathbb{H}^n$, cf.\ la définition~\ref{defi:point-conique}
plus bas).

De la même manière :
\begin{prop}
  \label{prop:group-conv-cocomp-dilatante-Hnbar}
  Soit $\Gamma$ un sous-groupe discret de $\mathbf{O}(n,1)$. Le groupe
  $\Gamma$ est convexe-cocompact si et seulement si, pour tout $z\in
  \Lambda_\Gamma$, l'action de $\Gamma$ est dilatante (dans la variété à bord
  $\overline{  \mathbb{H}}{}^n$ munie d'une quelconque distance
riemannienne) au point $z$. 
\end{prop}
À nouveau, le sens direct nécessite le lien entre dilatation au point $z$ et
convergence \og transverse\fg{}/\og conique\fg{} vers $z$. Pour la réciproque,
l'action de $\Gamma$ sur $\overline{ \mathbb{H}}{}^n \setminus \Lambda_\Gamma$
est propre (c'est toujours le cas pour un groupe discret) et la dilatation en
$\Lambda_\Gamma$ permet de montrer que cette action est cocompacte.

\subsection{Deux propriétés}
\label{sec:deux-proprietes}

Pour finir cette partie, énonçons deux propriétés importantes des groupes
convexes-cocompacts.

Tout d'abord, le groupe $\Gamma$ est quasi-isométrique au convexe $C\subset
\mathbb{H}^n$ et est donc hyperbolique au sens de Gromov. Il est possible de
réexprimer certains résultats de cette section, en y ajoutant une version
quantitative de la dilatation évoquée plus haut, sous une forme similaire au
théorème~\ref{th1}.

\begin{coro}
  \label{coro:group-conv-cocomp-hyp-dil-qi}
  Soit $\Gamma$ un sous-groupe de $\mathbf{O}(n,1)$. Le groupe $\Gamma$ est
  convexe-cocompact si et seulement si $\Gamma$ est hyperbolique au sens de
  Gromov et qu'il existe $\beta\colon \partial_\infty \Gamma \to \Lambda_\Gamma$ un
  homéomorphisme équivariant et 
  $k>0$, $K\geq 0$ tels que,
  pour tout rayon géodésique $(\gamma_p)_{p\in \mathbf{N}}$ dans $\Gamma$ avec
  $\gamma_0 = e_\Gamma$, on ait $\log \bigr\| \gamma_{p}^{-1} |_{
    T_{\beta(\gamma_\infty)} \partial_\infty \mathbb{H}^n}\bigr\| \geq k\, p-K
  $ où $\gamma_\infty \in \partial_\infty \Gamma$ est l'extrémité de la
  géodésique $(\gamma_p)_{p\in \mathbf{N}}$.
\end{coro}

Enfin les groupes convexes-cocompacts sont stables par petite déformation.
\begin{theo}
  \label{theo:group-conv-cocomp-stable-def}
  Soit $\Gamma$ un sous-groupe convexe-cocompact de $\mathbf{O}(n,1)$ et
  notons 
  $\iota\colon \Gamma \to \mathbf{O}(n,1)$ l'injection de $\Gamma$ dans
  $\mathbf{O}(n,1)$.

  Il existe alors un voisinage $U$ de $\iota$ dans l'ensemble
  $\mathrm{Hom}(\Gamma, \mathbf{O}(n,1))$ des homomorphismes de $\Gamma$ dans
  $\mathbf{O}(n,1)$ tel que, pour tout $\rho$ dans $U$, $\rho$ est injectif et
  son image $\rho(\Gamma)$ est convexe-cocompacte.
\end{theo}
La topologie sous-entendue sur $\mathrm{Hom}(\Gamma, \mathbf{O}(n,1))$ est
celle de la convergence simple. On peut en fait préciser les conclusions du
théorème : les actions de $\rho(\Gamma)$ et de $\Gamma = \iota(\Gamma)$ sur~$\overline{ \mathbb{H}}{}^n$ sont topologiquement conjuguées, l'ensemble
limite $\Lambda_{\rho(\Gamma)}$ varie continûment avec~$\rho$, etc.

\section{La préhistoire : absence de sous-groupes convexes- cocompacts en rang
  supérieur}
\label{sec:la-prehistoire-}

Les premiers résultats concernant les généralisations des groupes convexes-cocompacts en rang supérieur sont en fait négatifs : si l'on généralise
hâtivement la définition, aucun exemple intéressant et nouveau n'est
produit. Donnons ces résultats de rigidité pour les sous-groupes de
$\mathbf{SL}_n(\mathbf{R})$, $n\geq 3$.

Soit $X_n \subset \mathbf{M}_n(\mathbf{R})$ l'ensemble des matrices
symétriques, positives et de déterminant~$1$. C'est un espace homogène pour
l'action de $\mathbf{SL}_n(\mathbf{R})$ par transconjugaison : $(g,x) \mapsto
g\cdot x = g \, x \, {}^t\! g$ et les stabilisateurs sont compacts (celui de la
matrice identité $1_n\in X_n$ est le sous-groupe $\mathbf{SO}(n)$). Ainsi
$X_n$ admet une métrique riemannienne
$\mathbf{SL}_n(\mathbf{R})$-invariante. La courbure sectionnelle de cette
métrique (convenablement normalisée) varie entre $-1$ et $0$. Il y a donc dans
$X_n$ des directions \og hyperboliques\fg{} où la géométrie est semblable à
celle de l'espace hyperbolique mais aussi des directions \og plates\fg{} où la
géométrie est celle de l'espace euclidien. En particulier, il n'est pas
difficile de construire des quasi-géodésiques qui ne sont contenues dans le
voisinage métrique d'aucune géodésique (absence de lemme de Morse en rang
supérieur, voir cependant la partie~\ref{sec:hyperb-du-groupe}).
Le fait que~$X_n$ soit de courbure négative
implique qu'il est uniquement géodésique : toute paire de points est les
extrémités d'un unique segment géodésique. Une notion de convexité s'en déduit
aisément.

\begin{theo}[\cite{kleiner_leeb}]
  \label{theo:pas-de-cc-rang-sup}
  Soit $\Gamma$ un sous-groupe Zariski dense de $\mathbf{SL}_n(\mathbf{R})$,
  $n\geq 3$. Supposons qu'il existe un convexe $C\subset X_n$ non vide,
  $\Gamma$-invariant et sur lequel l'action de $\Gamma$ est propre et
  cocompacte.

  Alors $\Gamma$ est un réseau cocompact de $\mathbf{SL}_n(\mathbf{R})$.
\end{theo}
L'hypothèse de Zariski densité a été ajoutée ici pour simplifier la conclusion
mais il y a un résultat de structure sans cette hypothèse. Ce théorème vaut,
bien sûr, pour tout groupe de Lie semi-simple.

\smallskip

Pour énoncer le second résultat de rigidité, introduisons $D$ le sous-groupe
de $\mathbf{SL}_n(\mathbf{R})$ des matrices diagonales. Pour $\Gamma$ un
sous-groupe discret de $\mathbf{SL}_n(\mathbf{R})$, notons $E_\Gamma \subset \Gamma
\backslash \mathbf{SL}_n(\mathbf{R})$ l'adhérence des points fixes des
éléments réguliers (i.e.\ à valeurs propres deux à deux distinctes) de $D$
pour l'action à droite. C'est un ensemble $D$-invariant. Lorsque $n=2$, 
$\Gamma \backslash \mathbf{SL}_2(\mathbf{R})$ s'identifie au (à un revêtement
double du)  fibré unitaire tangent de $\Gamma \backslash \mathbb{H}^2$ et
l'action à droite de $D$ est le flot géodésique ; dans ce cas $E_\Gamma$ est
l'adhérence de la réunion des géodésiques fermées.
\begin{theo}[\cite{quint_cc}]
  \label{theo:pas-de-cc-rang-sup-flot}
  Soit $\Gamma$ un sous-groupe discret et Zariski dense de
  $\mathbf{SL}_n(\mathbf{R})$, $n\geq 3$. Supposons $E_\Gamma$ compact.
  Le groupe $\Gamma$ est alors un réseau cocompact de $\mathbf{SL}_n(\mathbf{R})$.
\end{theo}
Bien sûr, il y a, comme plus haut, la conclusion de structure sans l'hypothèse
de Zariski densité. Aussi le cas de sous-groupes discrets de groupes
algébriques sur les corps ultramétriques est traité. \textcite[\S
5]{quint_cc} déduit le théorème~\ref{theo:pas-de-cc-rang-sup} à partir de son
résultat.

\smallskip

À l'opposé, si l'on impose seulement au sous-groupe d'être de type fini et
quasi-isométriquement plongé, la classe des sous-groupes obtenus a de \og
mauvaises\fg{} propriétés :
\begin{itemize}
\item il existe un sous-groupe $\Gamma$ de $\mathbf{SL}_2(\mathbf{R}) \times
  \mathbf{SL}_2(\mathbf{R})$, isomorphe au groupe libre à deux générateurs,
  quasi-isométriquement plongé et vérifiant la propriété suivante : dans tout
  voisinage $U\subset \mathrm{Hom}( 
  \Gamma, \mathbf{SL}_2(\mathbf{R}) \times \mathbf{SL}_2(\mathbf{R}))$ de
  l'injection de $\Gamma$ dans $\mathbf{SL}_2(\mathbf{R}) \times
  \mathbf{SL}_2(\mathbf{R})$, il existe une représentation d'image dense
  \parencite[Prop.~A.1]{ggkw_anosov}.
\item il existe un sous-groupe de $\mathbf{SL}_2(\mathbf{R}) \times
  \mathbf{SL}_2(\mathbf{R})$, de type fini, quasi-isométriquement plongé mais
  qui n'est pas de présentation finie \parencite[Ex.~6.34]{kapovichleebporti14}.
\end{itemize}

\section{Les représentations Anosov}
\label{sec:les-repr-anos}

Cette partie introduit la notion, due à 
\textcite{labourie_anosov}, de
représentation Anosov. La définition ci-dessous est un peu éloignée de celle
donnée initialement par Labourie, elle correspond à la notion de sous-groupe
\og asymptotiquement plongé\fg{} de \textcite{kapovichleebporti14}. 
En particulier, l'absence de flot et de
décomposition hyperbolique pour ce flot dans les lignes à venir ne permet pas
d'apprécier le choix de la terminologie.

Fixons ici $\Gamma$ un groupe hyperbolique au sens de Gromov, son bord de
Gromov est toujours noté $\partial_\infty \Gamma$.

\begin{defi} 
  \label{defi:repr-anos}
  Une représentation $\rho\colon \Gamma \to \mathbf{GL}_n(\mathbf{R})$ est
  dite \emph{Anosov} s'il existe des applications continues et
  $\rho$-équivariante $\beta_1\colon \partial_\infty \Gamma \to
  \mathbb{P}^{n-1}(\mathbf{R})$ et $\beta_{n-1}\colon \partial_\infty \Gamma
  \to \mathbb{P}^{n-1*}(\mathbf{R})$ avec, pour tout $t\in \partial_\infty
  \Gamma$, $\beta_1(t)\!\subset \!\beta_{n-1}(t)$ et, pour tous 
  \mbox{$t\neq t' \in \partial_\infty \Gamma$, $\beta_1(t)\oplus \beta_{n-1}(t')\!=\!\mathbf{R}^n$}
  et vérifiant
  \begin{description}
  \item[Contraction] pour tous $t_1$ et $t_2$ dans $\partial_\infty \Gamma$,
    toute suite $(\gamma_p)_{p\in \mathbf{N}}$ de $\Gamma$ et tout \mbox{$\ell\in
    \mathbb{P}^{n-1}(\mathbf{R})$}, si $\lim \gamma_p = t_1$, $\lim
    \gamma_{p}^{-1} = t_2$ (limites dans $\Gamma \cup \partial_\infty \Gamma$)
    et si $\ell$ est transverse à $\beta_{n-1}(t_2)$, alors $\lim
    \rho(\gamma_p)\cdot \ell = \beta_1(t_1)$.
  \end{description}

  On dit qu'un sous-groupe $\Gamma$ de $\mathbf{GL}_n(\mathbf{R})$ est
  \emph{Anosov} si l'injection de $\Gamma$ dans $\mathbf{GL}_n(\mathbf{R})$
  est Anosov.
\end{defi}
Comme pour l'énoncé du théorème~\ref{th1}, seule l'image de $\Gamma$ dans
$\mathbf{PGL}_n(\mathbf{R})$ importe et l'on préfère considérer des
sous-groupes de $\mathbf{PGL}_n(\mathbf{R})$ ou de
$\mathbf{SL}_n(\mathbf{R})$.

La deuxième conclusion du théorème~\ref{th1} pourrait se réécrire \og
l'injection de $\Gamma$ dans $\mathbf{GL}_n(\mathbf{R})$ est Anosov\fg{}. La
propriété dynamique des suites divergentes donnée dans la définition est une
propriété de contraction faible. Celle donnée dans le théorème~\ref{th1} est
une propriété de contraction exponentielle. Dans la
partie~\ref{sec:diff-caract-des} seront données
plusieurs caractérisations des représentations Anosov, certaines ne faisant
pas a priori l'hypothèse d'hyperbolicité du groupe $\Gamma$. L'image finale
sera celle d'une classe de sous-groupes discrets généralisant la classe des
sous-groupes convexes-cocompacts.

Citons tout de suite trois propriétés des représentations Anosov qui
établissent déjà l'analogie avec les sous-groupes convexes-cocompacts :
\begin{itemize}
\item un sous-groupe Anosov est quasi-isométriquement plongé~;
\item l'action d'un sous-groupe Anosov $\rho(\Gamma)$ est dilatante sur
  $\mathbb{P}^{n-1}(\mathbf{R})$ en tout point de $\beta_1( \partial_\infty
  \Gamma)$~;
\item le sous-ensemble de $\mathrm{Hom}(\Gamma, \mathbf{GL}_n(\mathbf{R}))$
  constitué des représentations Anosov est ouvert.
\end{itemize}

La définition ci-dessus est un cas restreint de la notion. Il est possible de
faire la \og même\fg{} définition pour des représentations de $\Gamma$ dans
$G$ un groupe de Lie réductif ; la paire $( \mathbb{P}^{n-1}(\mathbf{R}),
\mathbb{P}^{n-1*}(\mathbf{R}))$ devant être remplacée par
$(G/P, G/P^{\mathrm{opp}})$ où $P$ est un sous-groupe parabolique de $G$ et
$P^{\mathrm{opp}}$ est un sous-groupe parabolique de $G$ opposé à $P$. On
utilisera parfois cette notion plus générale dans la suite et on
appellera les représentations la vérifiant $P$-Anosov ou $\mathcal{F}$-Anosov
avec $\mathcal{F}=G/P$. Plus encore, les caractérisations données plus bas
permettent d'englober les groupes d'isométries d'immeubles euclidiens.

\smallskip

\textcite{labourie_anosov} appelle plutôt ces représentations
$G/L$-Anosov, où $L=P\cap P^{ \mathrm{opp}}$ est le facteur de Levi, et les a
introduites pour les représentations des groupes de revêtements galoisiens de
variétés compactes munies d'un flot d'Anosov. L'idée de remplacer ces variétés
par le flot géodésique d'un groupe hyperbolique, et donc d'inclure tous les
groupes hyperboliques dans la définition, remonte à
\textcite{guichard_wienhard_dod}.

\section{Exemples}
\label{sec:exemples}

L'un des intérêts de la notion de représentation Anosov est l'abondance
d'exemples.

\subsection{Les espaces de Teichmüller généralisés}
\label{sec:les-espaces-de}

Il est aujourd'hui habituel d'appeler ainsi certaines composantes connexes de
la variété des représentations $\mathrm{Hom}(\Gamma_g, G)$ du groupe
fondamental $\Gamma_g \simeq \langle a_1, \dots, a_g, b_1, \dots, b_g \mid
[a_1, b_1]\cdots [a_g,b_g]\rangle$ d'une surface orientée, connexe, compacte,
sans bord et de genre $g\geq 2$, qui généralisent la composante de Teichmüller
qui est la composante connexe de $\mathrm{Hom}(\Gamma_g, \mathbf{SO}(2,1))$
constituée des représentations $\rho\colon \Gamma_g \to \mathbf{SO}(2,1)$
injectives, préservant l'orientation, 
et dont l'image $\rho(\Gamma_g)$ agit proprement sur $\mathbb{H}^2$ (on dira
que $\rho$ est une représentation fuchsienne) (une telle action est
automatiquement cocompacte : la surface quotient $\rho(\Gamma_g) \backslash
\mathbb{H}^2$ a son premier groupe d'homotopie isomorphe à $\Gamma_g$ et doit
donc être compacte). Il existe des versions des résultats ci-dessous pour des
surfaces ouvertes et leurs groupes fondamentaux (qui sont alors des groupes
libres de type fini) ; la discussion suivante se restreint, pour des raisons
de simplicité, aux surfaces fermées.

\subsubsection{Composante de Hitchin}
\label{sec:comp-de-hitch}

Soit $\tau_n\colon \mathbf{PGL}_2(\mathbf{R}) \simeq \mathbf{O}(2,1) \to
\mathbf{PGL}_n(\mathbf{R}) $ le morphisme induit par l'action de
$\mathbf{GL}_2( \mathbf{R})$ sur la puissance symétrique $n$-ième $S^{n-1}
\mathbf{R}^2 \simeq \mathbf{R}^n$ (explicitement $S^{n-1} \mathbf{R}^2$ est
l'espace des polynômes homogènes de degré $n-1$ en deux variables). La
composante de Hitchin est définie comme la composante connexe de
$\mathrm{Hom}(\Gamma_g, \mathbf{PGL}_n(\mathbf{R}))$ contenant les
compositions $\tau_n \circ \rho$ où $\rho\colon \Gamma_g \to
\mathbf{PGL}_2(\mathbf{R})$ est fuchsienne.

\begin{theo}[\cite{hitchin}]
  \label{theo:comp-de-hitch}
  Pour l'action de $\mathbf{PGL}_n(\mathbf{R})$ par conjugaison sur les
  représentations, la composante de Hitchin est
  $\mathbf{PGL}_n(\mathbf{R})$-difféomorphe à $\mathbf{R}^{ (2g-2)\times
    (n^2-1)} \times \mathbf{PGL}_n(\mathbf{R})$.
\end{theo}

\begin{theo}[\cite{labourie_anosov}]
  \label{theo:comp-de-hitch-anosov}
  Soit $B$ le sous-groupe de $\mathbf{PGL}_n(\mathbf{R})$ des matrices
  triangulaires supérieures. Toute représentation dans la composante de
  Hitchin est $B$-Anosov.
\end{theo}
Les travaux de Labourie apportent des informations plus précises sur
l'application au bord $\beta\colon \partial_\infty \Gamma_g \to
\mathbf{PGL}_n(\mathbf{R})/B \simeq \mathcal{F}lag(\mathbf{R}^n) = \{ (E_1,
\dots, E_{n-1}) \in \prod_{i=1}^{n} \mathrm{Gr}_i( \mathbf{R}^n) \mid \forall
i=1, \dots, n-2, \ E_i \subset E_{i+1}\}$ et permettent de donner une
\emph{caractérisation} des représentations qui appartiennent à la composante
de Hitchin. Les travaux de 
\textcite{fock_goncharov} donnent
également une caractérisation de ces représentations basée sur la notion de
positivité \og totale\fg{} des matrices (voir 
\cite{lusztigposred}).

{ Pour un groupe de Lie $G$ simple, adjoint et déployé sur $\mathbf{R}$, i.e.\
$G$ est égal à $\mathbf{PSL}_n(\mathbf{R})$, $\mathbf{PSp}_{2m}(\mathbf{R})$,
$\mathbf{PSO}(p,p+1)$ ou l'un des 5 groupes exceptionnels, le
$\mathbf{PSL}_2(\mathbf{R})$ principal de $G$ permet de définir la composante
de $\mathrm{Hom}(\Gamma_g, G)$ contenant les représentations fuchsiennes et
Hitchin a démontré que cette composante est $G$-difféomorphe à $\mathbf{R}^{
  (2g-2) \times \dim G} \times G$. En s'appuyant sur la notion de positivité
des représentations de $\Gamma_g$ dans $G$ développée par Fock et Goncharov, il
vient que toute représentation $\Gamma_g \to G$ dans la composante de Hitchin
est $B$-Anosov où $B$ est le sous-groupe de Borel de $G$.\par}

\subsubsection{Représentations maximales}
\label{sec:repr-maxim}

Le groupe fondamental de $\mathbf{Sp}_{2m}(\mathbf{R}) =\{ g \in
\mathbf{GL}_{2m}( \mathbf{R}) \mid {}^t\! g J_m g = J_m\}$, où $J_m = (
\begin{smallmatrix}
  0 & -1_m \\ 1_m & 0
\end{smallmatrix}
)$, est isomorphe à $\mathbf{Z}$. Le revêtement universel
$\widetilde{\mathbf{Sp}}_{2m}( \mathbf{R})$ de $\mathbf{Sp}_{2m}( \mathbf{R})$
est un groupe de Lie et le noyau de la projection naturelle
$\widetilde{\mathbf{Sp}}_{2m}( \mathbf{R}) \to \mathbf{Sp}_{2m}( \mathbf{R})$
est contenu dans le centre de $\widetilde{\mathbf{Sp}}_{2m}( \mathbf{R})$ et
s'identifie à $\pi_1( \mathbf{Sp}_{2m}(\mathbf{R})) \simeq
\mathbf{Z}$. 
Pour $g$ et
$h$ dans $\mathbf{Sp}_{2m}(\mathbf{R})$, le commutateur $\tilde{[}g,h\tilde{]}
= \tilde{g}\tilde{h}\tilde{g}{}^{-1}\tilde{h}{}^{-1}$ de relevés $\tilde{g}$
et $\tilde{h}\in \widetilde{\mathbf{Sp}}_{2m}( \mathbf{R})$ de $g$ et $h$ ne
dépend pas de ces relevés et sa projection dans $\mathbf{Sp}_{2m}(
\mathbf{R})$ est égale à $[g,h]$.

Soit maintenant $\rho\colon \Gamma_g \to \mathbf{Sp}_{2m}(\mathbf{R})$ un
morphisme, la projection de $\prod_{i=1}^{g} \tilde{[} \rho(a_i), \rho(b_i)
\tilde{]}$ dans $\mathbf{Sp}_{2m}( \mathbf{R})$ est égale à $\prod_{i=1}^{g} [
\rho(a_i), \rho(b_i)] = \rho\bigl( \prod_{i=1}^{g} [ a_i, b_i]\bigr) = \rho(
e_\Gamma) = e_{\mathbf{Sp}_{2m}( \mathbf{R})}$ et cet élément appartient donc
au noyau de $\widetilde{\mathbf{Sp}}_{2m}( \mathbf{R}) \to \mathbf{Sp}_{2m}(
\mathbf{R})$, c'est-à-dire à $\pi_1( \mathbf{Sp}_{2m}(\mathbf{R})) \simeq
\mathbf{Z}$.
\begin{defi}
  \label{defi:nombre-euler}
  On appelle \emph{nombre d'Euler} ou \emph{nombre de Toledo} et on note
  $\mathbf{e}(\rho)$ le nombre entier $\prod_{i=1}^{g} \tilde{[} \rho(a_i),
  \rho(b_i) \tilde{]}$.
\end{defi}
\begin{theo}[Inégalité de Milnor et Wood]
  \label{theo:milnor-wood} 
  Pour toute représentation $\rho\colon \Gamma_g \to \mathbf{Sp}_{2m}(
  \mathbf{R})$, on a $| \mathbf{e}(\rho)| \leq m(g-1)$.
\end{theo}
\begin{defi}
  \label{defi:repr-maxim}
  Une  représentation $\rho\colon \Gamma_g \to \mathbf{Sp}_{2m}(
  \mathbf{R})$ est dite \emph{maximale} si $ |\mathbf{e}(\rho)| = m(g-1)$.
\end{defi}
Soit $Q$ le stabilisateur dans
$\mathbf{Sp}_{2m}( \mathbf{R})$ de l'espace engendré par $ e_1, \dots, e_m$ (toujours avec $(e_i)_{i=1, \dots, 2m}$ la
base canonique de $\mathbf{R}^{2m}$), c'est un sous-groupe parabolique
maximal de $\mathbf{Sp}_{2m}( \mathbf{R})$. 
\begin{theo}
  [\cite{burger_iozzi_labourie_wienhard}]
  Toute représentation maximale est $Q$-Anosov.
\end{theo}
Ce résultat est basé sur les travaux de 
\textcite{burger_iozzi_wienhard_toledo} établissant l'existence d'une application
équivariante, continue à droite $\partial_\infty \Gamma_g \to
\mathbf{Sp}_{2m}( \mathbf{R})/ Q \simeq \mathcal{L}ag( \mathbf{R}^{2m}) = \{
h\in \mathrm{Gr}_m( \mathbf{R}^{2m}) \mid \forall x,y \in h,\ {}^t\! x J_m y=0
\}$ qui satisfait de surcroît une propriété de \og maximalité\fg{}.

{ Lorsque $G$ est un groupe de Lie simple et de type hermitien, le groupe
fondamental de $G$ a toujours un facteur cyclique infini ce qui permet de
définir un nombre d'Euler/Toledo pour les représentations $\Gamma_g \to G$.
Il y a encore une majoration de cet entier et donc une notion de
représentations maximales. Les travaux de Burger, Iozzi et Wienhard montrent
que ces représentations sont Anosov par rapport à un sous-groupe parabolique
bien déterminé.\par
}
\subsection{Autres exemples}
\label{sec:autres-exemples}

\subsubsection{Les sous-groupes convexes-cocompacts}
\label{sec:les-groupes-convexes}

Ces sous-groupes sont Anosov ! Plus généralement si $\Gamma$ est un sous-groupe
convexe-cocompact de $\mathbf{O}(n,1)$ et si \mbox{$\tau\colon \mathbf{O}(n,1)\! \to \!G$} est
un plongement, alors $\tau(\Gamma)$ est Anosov relativement à un sous-groupe
parabolique de $G$ déterminé par le morphisme $\tau$.

\subsubsection{Groupes de Schottky}
\label{sec:groupes-de-schottky}

Les groupes de Schottky, c'est à dire les sous-groupes engendrés par une
paire (de puissances) d'éléments réguliers et suffisamment transverses, sont Anosov. Une démonstration
géométrique, i.e.\ à partir de l'action sur l'espace symétrique et n'utilisant
pas de versions du lemme du ping-pong, se trouve
dans \cite{kapovichleebporti14} et est basée sur le
théorème~\ref{theo:local-global}. 

\subsubsection{Groupes de surfaces dans $\mathbf{PGL}_3(\mathbf{R})$}
\label{sec:groupes-de-surfaces}

\textcite{barbot_anosov} 
étudie les représentations Anosov $\rho$ de groupes
de surface $\Gamma_g$ dans $\mathbf{PGL}_3( \mathbf{R})$. Il montre en
particulier qu'il y a toujours un ouvert $\Omega \subset \mathcal{F}lag( \mathbf{R}^3)$ de
la variété drapeau sur lequel $\rho(\Gamma_g)$ agit proprement et avec
quotient compact. Il caractérise aussi les représentations de la composante de
Hitchin avec $\Omega$ : $\rho$ appartient à la composante de Hitchin si et
seulement si $\Omega$ n'est pas connexe.

\subsubsection{Convexes divisibles}
\label{sec:convexes-divisibles}

On appelle ainsi les ouverts proprement convexes $\Omega$ de l'espace projectif
$\mathbb{P}^{n-1}( \mathbf{R})$ pour lesquels il existe un sous-groupe
$\Gamma$ de $\mathrm{Aut}(\Omega) =\{ g \in \mathbf{PGL}_n( \mathbf{R}) \mid
g\cdot \Omega = \Omega \}$ agissant proprement sur $\Omega$ avec quotient
compact. Les réseaux de $\mathbf{O}(n-1,1)$ sont de tels exemples avec $\Omega
= \mathbb{H}^{n-1}$. 
\textcite{benoist_cd1} démontre que le groupe
$\Gamma$ est hyperbolique au sens de Gromov si et seulement si $\Omega$ est
strictement convexe, auquel cas le bord de Gromov $\partial_\infty \Gamma$
s'identifie de manière $\Gamma$-équivariante à~$\partial \Omega$. En utilisant
l'action sur le convexe polaire $\Omega^0 \subset \mathbb{P}^{n-1*}(
\mathbf{R})$ on obtient une seconde application $\partial_\infty \Gamma \to
\mathbb{P}^{n-1*}( \mathbf{R})$. À partir de là, il est aisé de
démontrer que $\Gamma \to \mathbf{PGL}_n( \mathbf{R})$ est Anosov.

Le théorème de stabilité des représentations Anosov peut être utilisé dans ce
cas pour retrouver un résultat de 
\textcite{koszul} de stabilité des
structures projectives convexes sur une variété compacte.

\textcite{kapovich_convex} 
construit des ouverts convexes divisibles
$( \Omega, \Gamma)$ pour lesquels le quotient $\Gamma \backslash \Omega$ est
une variété de Gromov et Thurston ; ceci donne les premiers exemples de
sous-groupes Anosov qui ne sont pas isomorphes à un réseau d'un groupe de Lie.

\subsubsection{Sous-groupes quasi-fuchsiens dans $\mathbf{SO}(2,n)$}
\label{sec:sous-groupes-quasi}
\textcite{barbot_merigot_fusionpubli} 
ont défini une notion de sous-groupe $\Gamma$ quasi-fuchsien
dans $\mathbf{SO}(2,n)$ : l'ensemble limite~$\Lambda_\Gamma$, dans le  projectivisé
du cône isotrope pour la forme quadratique $q_{2,n}$, est homéomorphe à la
sphère de dimension $n-1$ et vérifie une propriété d'\og acausalité\fg{}. Ils démontrent 
qu'alors $\Gamma$ est Anosov et que,
réciproquement, 
un sous-groupe Anosov dont le bord de
Gromov $\partial_\infty \Gamma$ est homéomorphe à une sphère de dimension
$n-1$ est quasi-fuchsien. 
\textcite{barbotdefadsqf} démontre ensuite que le sous-ensemble
des représentations quasi-fuchsiennes dans $\mathrm{Hom}( \Gamma, \mathbf{SO}(
2,n))$ est fermé ; comme cet ensemble est également ouvert, il est réunion de composantes connexes.

\subsubsection{Groupes de Coxeter}
\label{sec:groupes-de-coxeter}

En se basant sur leurs travaux qui développent la notion de sous-groupes
convexes-cocompacts en géométrie projective (voir le paragraphe~\ref{sec:le-retour-des} plus bas), 
\textcite{danciger_gueritaud_kassel_psdo-hyp,danciger_gueritaud_kassel-cc-proj} montrent que tous
les groupes de Coxeter à angles droits et hyperboliques admettent des
représentations Anosov ; les images de ces représentations sont des groupes de
réflexions hyperplanes. 
Ceci permet d'avoir de nouveaux exemples de représentations
Anosov pour des groupes discrets qui ne sont pas des réseaux d'un groupe de Lie.

En utilisant ces travaux, 
\textcite{lee_marquis} ont
donné les premiers exemples de sous-groupes quasi-fuchsiens (au sens de Barbot
et Mérigot, voir \S~\ref{sec:sous-groupes-quasi}) qui ne sont pas isomorphes à des réseaux de $\mathbf{O}(1,n)$.

\section{Différentes caractérisations des sous-groupes Anosov}
\label{sec:diff-caract-des}

Cette partie passe en revue quelques caractérisations des représentations
Anosov (ici des représentations $P_1$-Anosov dans $\mathbf{GL}_n(
\mathbf{R})$), il y en a bien d'autres et le lecteur pourra consulter
\textcite{kapovichleeb15,kapovichleeb17,kapovichleebporti14,kapovichleebporti14_2,kapovichleebporti15,kapovichleebportiv2,kapovichleebporti17,ggkw_anosov,guichard_wienhard_dod,labourie_anosov}
pour un tableau plus complet.

\subsection{Un peu de géométrie de l'espace symétrique}
\label{sec:un-peu-de}

Ce paragraphe aborde de manière succincte, pour le groupe $\mathbf{SL}_n( \mathbf{R})$, la géométrie de l'espace symétrique
et les variétés drapeaux comme \og bord à l'infini\fg{} de l'espace symétrique
; ceci est traité de manière plus systématique dans les articles de Kapovich,
Leeb et Porti.

Nous travaillerons dans cette partie plutôt avec le groupe $\mathbf{SL}_n(
\mathbf{R})$ que $\mathbf{GL}_n( \mathbf{R})$ ; il est commode et naturel de
se \og débarrasser\fg{} du facteur central de $\mathbf{GL}_n( \mathbf{R})$ et
de travailler avec un groupe de Lie connexe.  L'espace symétrique
$\mathbf{SL}_n( \mathbf{R}) / \mathbf{SO}(n) $ s'identifie à $X_n =\{ s \in
\mathbf{M}_n( \mathbf{R}) \mid {}^t\! s =s, \, \det s=1, \text{ et } 
\forall x\in \mathbf{R}^n \setminus \{0\},\, {}^t\! x s x >0\}$. Soit $s_0 = 1_n$
le point de $X_n$ dont le stabilisateur est $\mathbf{SO}(n)$.

La décomposition
de Cartan dans $\mathbf{SL}_n( \mathbf{R})$ (qui est reliée aux valeurs
principales déjà mentionnées) implique le (et est même équivalente au)  fait
que pour tout $s$ dans $X_n$, il existe un élément $k$ de $\mathbf{SO}(n)$ et
une matrice diagonale ${\boldsymbol{\mu}}= \mathrm{diag}( \mu_1, \dots, \mu_n)$ avec
$\mu_1 \geq \mu_2 \geq \cdots \geq \mu_n$ et $\sum_{i=1}^{n} \mu_i =0$ tels
que $s = k e^{{\boldsymbol{\mu}}} \cdot s_0$. Le $n$-uplet $(\mu_1, \dots, \mu_n)$
est alors uniquement déterminé. On dira que le segment $s_0 s$ est
\emph{régulier} si à la fois $\mu_1> \mu_2$ et $\mu_{n-1}> \mu_n$. Plus
généralement, si $s$ et $s'$ sont deux points de $X_n$, il existe un élément
$g$ de $\mathbf{SL}_n( \mathbf{R})$ et une matrice diagonale ${\boldsymbol{\mu}} =
\mathrm{diag}( \mu_1, \dots, \mu_n)$ avec $\mu_1 \geq \cdots \geq \mu_n$ tels
que $s= g\cdot s_0$ et $s' = g e^{{\boldsymbol{\mu}}} \cdot s_0$. La matrice
${\boldsymbol{\mu}}$ est uniquement déterminée par $(s,s')$. Le segment $ss'$ est
dit \emph{régulier} si $\mu_1 > \mu_2$ et $\mu_{n-1}> \mu_n$. La matrice
associée à $(s',s)$ est la matrice diagonale $\mathrm{diag}(-\mu_n,
-\mu_{n-1}, \dots, -\mu_1)$. Ainsi le segment $s's$ est régulier si et
seulement si $ss'$ l'est.

La variété des drapeaux constitués d'une droite et d'un hyperplan est notée
$\mathcal{F}_{1,n-1} = \{ (\ell,h) \in \mathbb{P}^{n-1}( \mathbf{R}) \times
\mathbb{P}^{n-1*}( \mathbf{R}) \mid \ell\subset h\}$. C'est un espace homogène
sous l'action de $\mathbf{SO}(n)$ et le stabilisateur dans $\mathbf{SO}(n)$ de
$f_0 = (\ell_0, h_0)$ où $\ell_0 = \mathbf{R} e_1$ et $h_0 =\mathbf{R} e_1
\oplus \cdots \oplus \mathbf{R}
e_{n-1}$ est le sous-groupe $M$ de $\mathbf{SO}(n)$ constitué des
matrices diagonales par blocs, $M = \Bigl\{ \Bigl(
\begin{smallmatrix}
  \varepsilon & & \\ & a & \\ & & \delta
\end{smallmatrix}
\Bigr) \ \Big|\ \varepsilon, \delta \in \{ \pm 1\},\, a\in \mathbf{O}(n-2),\,
\varepsilon \delta \det a =1\Bigr\}$.

On appellera \emph{cône de Weyl} (centré en $s_0$ et de direction $f_0$) et on
notera $V(s_0, f_0)$ le sous-ensemble de $X_n$ formé des $m e^{{\boldsymbol{\mu}}}
\cdot s_0$ avec ${\boldsymbol{\mu}}= \mathrm{diag}( \mu_1, \dots, \mu_n)$,
$\mu_1 \geq \mu_2 \geq \cdots \geq \mu_n$ et $m\in M$. Si $f=(\ell,h)$ appartient
à $\mathcal{F}_{1,n-1}$, soit $k\in \mathbf{SO}(n)$ tel que $f = k\cdot f_0$,
le \emph{cône de Weyl} $V(s_0, f)$ est l'image de $V(s_0, f_0)$ par $k$ :
$V(s_0, f) = k\cdot V(s_0, f_0) =\{ k_1 e^{{\boldsymbol{\mu}}} \cdot s_0 \mid
{\boldsymbol{\mu}}= \mathrm{diag}( \mu_1, \dots, \mu_n),\, \mu_1 \geq \cdots \geq
\mu_n,\, k_1\in \mathbf{SO}(n),\, k_1 \cdot f_0 = f\}$. Un point $s = g\cdot s_0$ appartient à
$V(s_0,f)$ si et seulement si la matrice $ g\, {}^t\! g$ stabilise $\ell$ et $h$
(et donc $\ell^\perp$ et $h^\perp$) et si $\| g\, {}^t\! g |_\ell \| \geq \| g\,
{}^t\! g |_{\ell^\perp} \|$ et $\| (g\, {}^t\! g)^{-1} |_{h^\perp} \| \geq \|
(g\, {}^t\! g)^{-1} |_{h} \|$. Il existe donc toujours un cône de Weyl $V(s_0,
f)$ ($f=(\ell,h)\in \mathcal{F}_{1,n-1}$) contenant $s= g\cdot s_0 = k
e^{{\boldsymbol{\mu}}} \cdot s_0$ et ce cône de Weyl est uniquement déterminé à
$s$ si et seulement si $\| g\, {}^t\! g |_\ell \| > \| g\,
{}^t\! g |_{\ell^\perp} \|$ et $\| (g\, {}^t\! g)^{-1} |_{h^\perp} \| > \|
(g\, {}^t\! g)^{-1} |_{h} \|$. Comme ici $e^{\mu_1} = \| g\, {}^t\! g |_\ell \|$,
$e^{ \mu_2} =  \| g\, {}^t\! g |_{\ell^\perp} \|$, $e^{ -\mu_n} = \| (g\, {}^t\!
g)^{-1} |_{h^\perp} \|$ et $e^{-\mu_{n-1}} = \| (g\, {}^t\! g)^{-1} |_{h} \|$,
l'unicité de ce cône de Weyl est équivalente à ce que le segment $s_0 s$
est régulier.

On peut bien sûr considérer des cônes de Weyl de sommets quelconques : si
$s\in X_n$ et $f\in \mathcal{F}_{1,n-1}$, il existe $g\in \mathbf{SL}_n(
\mathbf{R})$ tel que $s=g\cdot s_0$ et $f=g\cdot f_0$ (toujours car le
stabilisateur de $s_0$ dans $\mathbf{SL}_n( \mathbf{R})$ agit transitivement
sur $\mathcal{F}_{1,n-1}$), on pose $V(s,f) = g\cdot V(s_0, f_0)$. Comme $g$
est uniquement déterminé par multiplication à droite par un élément de $M$
près et comme $V(s_0, f_0)$ est invariant par $M$, cet ensemble dépend
uniquement de la paire $(s,f)$  et non de $g$. Pour tout autre point $s'$ de
$X_n$, il existe un cône de Weyl $V(s,f)$ contenant $s'$ et ce cône est unique
si et seulement si $ss'$ est un segment régulier.

Soit maintenant $ss'$ un segment régulier, le \emph{diamant}
$\diamondsuit_{s,s'}$ déterminé par $ss'$ est l'intersection du cône de Weyl
$V(s,f)$ contenant $s'$ et du cône de Weyl $V(s',f')$ contenant~$s$ :
$\diamondsuit_{s,s'} = V(s,f) \cap V(s',f')$. C'est un sous-ensemble
\emph{compact} de $X_n$. On a $\diamondsuit_{s',s} = \diamondsuit_{s,s'}$ et,
pour tout $s'' \in \diamondsuit_{s,s'}$ tel que $ss''$ est régulier,
$\diamondsuit_{s,s''} \subset \diamondsuit_{s,s'}$. 
Notons aussi que, sous ces hypothèses, les drapeaux $f=(\ell,h)$ et $f'=(\ell',h')$
sont \emph{transverses} : $\ell\cap h'=0$ et $\ell' \cap h=0$.

Nous aurons également besoin de versions quantitatives de la régularité. Un
segment~$ss'$ de~$X_n$ sera dit \emph{$\epsilon$-régulier} (où $\epsilon>0$
est fixé) si $s=g\cdot s_0$, $s' = ge^{{\boldsymbol{\mu}}} \cdot s_0$ avec
${\boldsymbol{\mu}} = \mathrm{diag}(\mu_1, \dots, \mu_n)$, $\mu_1 \geq \cdots \geq
\mu_n$, sont tels que $\mu_1 -\mu_2 \geq \epsilon\, d_{X_n}(s,s')$ et $\mu_{n-1} -\mu_{n} \geq
\epsilon\, d_{X_n}(s,s')$ avec $d_{X_n}(s,s') = \sqrt{ \sum_{i=1}^{n} \mu_{i}^{2}}
$ ($d_{X_n}$ est la distance riemannienne dans
$X_n$). On notera $V^\epsilon(s,f)$ l'ensemble des $s'$ dans $V(s,f)$ tels que
$ss'$ est $\epsilon$-régulier. On notera $\diamondsuit_{s,s'}^{\epsilon}$
l'intersection $V^\epsilon(s,f) \cap V^\epsilon(s',f')$ où $f$ et $f'$ sont
déterminés par le segment régulier $ss'$ comme ci-dessus.

\subsection{Un peu de dynamique sur $\mathcal{F}_{1,n-1}$}
\label{sec:un-peu-de-dyn}

Les notions introduites dans le paragraphe précédent permettent de définir
l'ensemble limite dans $\mathcal{F}_{1,n-1}$ d'un sous-groupe $\Gamma$ de
$\mathbf{SL}_n( \mathbf{R})$.

\begin{defi}
  \label{defi:ens-lim}
  On note $\Lambda_\Gamma \subset \mathcal{F}_{1,n-1}$ et on appelle
  \emph{ensemble limite} de $\Gamma$ dans $\mathcal{F}_{1,n-1}$ l'ensemble des
  limites de suites $(f_p)_{p\in \mathbf{N}}$ de $\mathcal{F}_{1,n-1}$ telles
  qu'il existe une suite $(\gamma_p)_{p\in \mathbf{N}}$ dans $\Gamma$ avec
  \begin{itemize}
  \item $\lim \mu_1(\gamma_p) - \mu_2( \gamma_p) = \lim \mu_{n-1}(\gamma_p) -
    \mu_n( \gamma_p) = +\infty$,
  \item $\lim d_{X_n}( s_0, \gamma_p \cdot s_0) = +\infty$ et, 
  \item pour tout
    $p\in \mathbf{N}$, $\gamma_p \cdot s_0 \in V(s_0, f_p)$.
  \end{itemize}

\end{defi}
L'ensemble limite est fermé, il ne dépend pas du point base $s_0$. On
pourrait, dans la définition, supposer seulement que $\gamma_p \cdot s_0$ est
dans un $R$-voisinage de $V(s_0, f_p)$ où $R\geq 0$ est indépendant de $p$. On
appelle points coniques les points de $\Lambda_\Gamma$ obtenus ainsi en
prenant une suite $(f_p)_{p\in \mathbf{N}}$ constante :
\begin{defi}
  \label{defi:point-conique}
  Un élément $f$ de $\mathcal{F}_{1,n-1}$ est un \emph{point limite conique}
  s'il existe une suite $(\gamma_p)_{p\in \mathbf{N}}$ de $\Gamma$ et $R\geq
  0$ avec $\lim d_{X_n}( s_0, \gamma_p \cdot s_0) = +\infty$ et, pour tout
  $p\in \mathbf{N}$, $\gamma_p \cdot s_0$ appartient au $R$-voisinage du cône
  de Weyl $V(s_0, f)$.
\end{defi}

\smallskip

Des caractérisations des sous-groupes Anosov seront aussi données en termes de
l'action de $\Gamma$ sur $\mathcal{F}_{1,n-1}$, i.e.\ sans référence à
l'espace symétrique $X_n$. Les notions pertinentes sont les suivantes.

\begin{defi}
  \label{defi:contrac-convergence}
  Une suite $(g_p)_{p\in \mathbf{N}}$ de $\mathbf{SL}_n( \mathbf{R})^{
    \mathbf{N}}$ est dite \emph{contractante} sur $\mathcal{F}_{1,n-1}$ s'il
  existe $f^+$ et $f^-$ dans $\mathcal{F}_{1,n-1}$ tels que, pour tout $f\in
  \mathcal{F}_{1,n-1}$, si $f$ est transverse à $f^-$, alors la suite $(g_p
  \cdot f)_{p\in \mathbf{N}}$ tend vers $f^+$.

  L'action de $\Gamma$ sur $\mathcal{F}_{1,n-1}$ sera dite \emph{de
    convergence} si toute suite de $\Gamma$ tendant vers l'infini a une
  sous-suite contractante. (On dit qu'une suite $(\gamma_p)_{p\in \mathbf{N}}$
  tend vers l'infini quand $\lim d_{X_n}( s_0, \gamma_p \cdot s_0) =
  +\infty$.)
\end{defi}
Si $(g_p)_{p\in \mathbf{N}}$ est une suite contractante, alors les drapeaux
$f^\pm$ de la définition sont uniquement déterminés et la convergence vers
$f^+$ est en fait uniforme sur les compacts de $\mathcal{F}_{1,n-1}$ contenus
dans l'ouvert des éléments de $\mathcal{F}_{1,n-1}$ transverses à $f^-$.

On peut détecter l'existence d'une sous-suite contractante à l'aide de la
projection de Cartan. Pour tout $g$ dans $\mathbf{SL}_n(\mathbf{ R})$ notons
encore ${\boldsymbol{\mu}}(g) = \mathrm{diag}( \mu_1(g), \dots, \mu_n(g))$ l'unique
matrice diagonale à coefficients décroissants contenue dans la double classe
$\mathbf{SO}(n) g \mathbf{SO}(n)$. Alors une suite $( g_p)_{p\in \mathbf{N}}
\in \mathbf{SL}_n( \mathbf{R})^{ \mathbf{N}}$ a une sous-suite contractante si
et seulement si $\sup \{ \mu_1(g_p) -\mu_2( g_p), \, p\in \mathbf{N} \}= \sup
\{ \mu_{n-1}(g_p) -\mu_n( g_p), \, p\in \mathbf{N} \} = + \infty$.

\smallskip

Les définitions ci-dessus sont empruntées des articles de Kapovich, Leeb et
Porti. 
\textcite{benoist} définit aussi une notion d'ensemble limite
qu'il étudie en détails. Dans ses travaux sur les mesures de Patterson et
Sullivan en rang supérieur, 
\textcite{albuquerque} introduit également
une notion de point limite conique.

\subsection{Sous-groupes de convergence, dilatant et à ensemble limite
  transverse}
\label{sec:sous-groupes-de}

La caractérisation donnée ici des représentations Anosov est un analogue du
théorème~\ref{theo:group-conv-cocomp-dilatante-dHn} pour les sous-groupes convexes-cocompacts.

\begin{theo}[\cite{kapovichleebporti14}]
  \label{theo:trans-conv-dil}
  Soit $\Gamma$ un sous-groupe discret de $\mathbf{SL}_n( \mathbf{R})$. On suppose
  \begin{description}
  \item[Transversalité] Pour tous $f\neq f'$ dans $\Lambda_\Gamma$, les
    drapeaux $f$ et $f'$ sont transverses.
  \item[Convergence] L'action de $\Gamma$ sur $\mathcal{F}_{1,n-1}$ est de
    convergence (définition~\ref{defi:contrac-convergence}).
  \item[Dilatation] L'action de $\Gamma$ en tout point $f\in \Lambda_\Gamma$
    est dilatante : il existe $\gamma\in\Gamma$, $c>0$ et un voisinage $U$ de
    $f$ dans $\mathcal{F}_{1,n-1}$ tels que, pour tous $f_1$, $f_2\in U$, $d_{
    \mathcal{F}_{1,n-1}} ( \gamma\cdot f_1, \gamma\cdot f_2) \geq c\,d_{
    \mathcal{F}_{1,n-1}} ( f_1, f_2)$ ($d_{  \mathcal{F}_{1,n-1}}$ est une
  distance riemannienne sur $\mathcal{F}_{1,n-1}$).
  \end{description}
  Alors le groupe~$\Gamma$ est de type fini, hyperbolique au sens de Gromov et
  $\Gamma$ est un sous-groupe Anosov de $\mathbf{SL}_n( \mathbf{R})$. De plus,
  si $\beta_1\colon \partial_\infty \Gamma \to \mathbb{P}^{n-1}( \mathbf{R})$
  et $\beta_{n-1}\colon \partial_\infty \Gamma \to \mathbb{P}^{n-1*}(
  \mathbf{R})$ sont les applications au bord associées, l'application $\beta =
  (\beta_1, \beta_{n-1})\colon \partial_\infty \Gamma \to
  \mathcal{F}_{1,n-1}\subset  \mathbb{P}^{n-1}( \mathbf{R}) \times
  \mathbb{P}^{n-1 *}( \mathbf{R})$ induit un homéomorphisme
  $\Gamma$-équivariant sur $\Lambda_\Gamma$.
\end{theo}

Dans ce théorème, l'hyperbolicité est déduite du résultat suivant :

\begin{theo}[\cite{bowditch_tophyp}]
  \label{theo:conv-grp-hyp}
  Soit $Z$ un espace topologique compact, métrisable et parfait et soit
  $\Gamma$ un groupe d'homéomorphismes de $Z$ agissant proprement sur $\{
  (z_1, z_2, z_3) \in Z^3 \mid z_1 \neq z_2,\, z_2 \neq z_3,\, z_3 \neq z_1\}$
  et avec quotient compact. Le groupe~$\Gamma$ est alors de type fini,
  hyperbolique au sens de Gromov et son bord de Gromov $\partial_\infty
  \Gamma$ s'identifie à $Z$ de manière $\Gamma$-équivariante.
\end{theo}

\subsection{Sous-groupes réguliers à ensemble limite conique et transverse}
\label{sec:sous-group-regul}

Cette caractérisation entre en résonance avec le
corollaire~\ref{coro:group-conv-cocomp-hyp-dil-qi} pour les sous-groupes convexes-cocompacts.

\begin{theo}[\cite{kapovichleebporti14}]
  \label{theo:reg-coni-trans}
  Soit $\Gamma$ un sous-groupe de $\mathbf{SL}_n( \mathbf{R})$. On suppose
  \begin{description}
  \item[Régularité] $\lim_{\gamma \to \infty, \gamma\in\Gamma} \mu_1(\gamma) -
    \mu_2(\gamma) =+\infty$ (de manière équivalente, pour toute suite
    $(\gamma_p)_{p\in \mathbf{N}}$ d'éléments deux à deux distincts, $\lim_p
    \mu_1(\gamma_p) - \mu_2(\gamma_p) =+\infty$).
  \item[Transversalité] Pour tous $f$, $f'$ dans $\Lambda_\Gamma$, si $f\neq
    f'$, alors $f$ et $f'$ sont transverses.
  \item[Conicalité] Tout point $f$ de $\Lambda_\Gamma$ est un point limite
    conique (définition~\ref{defi:point-conique}).
  \end{description}
  Alors le groupe $\Gamma$ est de type fini, hyperbolique au sens de Gromov et
  l'injection de $\Gamma$ dans $\mathbf{SL}_n( \mathbf{R})$ est une
  représentation Anosov.
\end{theo}

Bien sûr, réciproquement, un sous-groupe Anosov satisfait à toutes les
conditions énoncées dans ces deux théorèmes.

\subsection{Action de type Morse}
\label{sec:action-de-type}

Le lemme de Morse, pour les espaces hyperboliques, affirme que les
quasi-géodésiques sont à distance finie de géodésiques. Un tel énoncé n'est pas
vrai en rang supérieur mais 
\textcite{kapovichleebporti14_2} en ont donné
une version pour $X_n$ (et pour les espaces symétriques et les immeubles
euclidiens) : un segment géodésique \emph{régulier} (i.e.\ vérifiant la
condition énoncée dans le théorème~\ref{th1}) est proche d'un diamant
(respectivement, un rayon géodésique régulier est proche d'un
 cône de Weyl et une géodésique biinfinie est proche de la réunion de deux
 cônes de Weyl de même sommet et de directions transverses). La
partie~\ref{sec:hyperb-du-groupe} 
suivante donne quelques énoncés intermédiaires qui
mènent à ce résultat d'approximation des quasi-géodésiques régulières.

Voyons maintenant qu'un sous-groupe hyperbolique vérifiant cette propriété
d'approximation est Anosov.

\begin{theo}[\cite{kapovichleebporti14}]
  \label{theo:morse-est-anosov}
  Soit $\Gamma$ un sous-groupe de type fini de $\mathbf{SL}_n(
  \mathbf{R})$. On suppose que $\Gamma$ est hyperbolique au sens de Gromov et que
  \begin{description}
  \item[Régularité]  $\lim_{\gamma \to \infty, \gamma\in\Gamma} \mu_1(\gamma) -
    \mu_2(\gamma) =+\infty$.
  \item[Morse] Il existe $R\!\geq\! 0$ pour lequel la propriété suivante est
    vérifiée : pour toute géodésique finie $(
    \gamma_p)_{p=0,\dots, P}$ dans le graphe de Cayley de $\Gamma$, il existe
    $s$ et $s'$ dans $X_n$ tels que $d_{X_n}( \gamma_0\cdot s_0, s)\!\leq \!R$,
    \mbox{$d_{X_n}( \gamma_P\cdot s_0, s')\!\leq \!R$}, le segment $ss'$ est
    régulier et, 
    pour tout $p=0,\dots, P$, $\gamma_p\cdot s_0$ appartient au \mbox{$R$-voisinage}
    du diamant $\diamondsuit_{s,s'}$.
   \end{description}
   Le groupe $\Gamma$ est alors un sous-groupe Anosov.
\end{theo}

Réciproquement, un sous-groupe Anosov est \og Morse\fg{} de manière
quantitative :
\begin{theo}[\cite{kapovichleebporti14}]
  \label{theo:anosov-est-morse}
  Soit $\Gamma$ un sous-groupe Anosov de $\mathbf{SL}_n( \mathbf{R})$. Il
  existe alors $L\geq 1$, $A\geq 0$, $R\geq 0$ et $\epsilon>0$ tels que, pour
  toute géodésique finie $(\gamma_p)_{p=0,\dots, P}$ dans le graphe de Cayley
  de $\Gamma$, on a :
  \begin{itemize}
  \item l'application $\{0,1,\dots,P\} \to X_n / p \mapsto \gamma_p\cdot s_0$
    est une $(L,A)$-quasi-géodésique ;
  \item il existe $s \in B_{X_n}( \gamma_0\cdot s_0, R)$ et $s' \in B_{X_n}(
    \gamma_P\cdot s_0, R)$ tels que $\{ \gamma_p\cdot s_0, \, p=0,\dots, P\}$
    est contenu dans le $R$-voisinage de $\diamondsuit^{\epsilon}_{s,s'}$.
  \end{itemize}
\end{theo}
Lorsque les conclusions du théorème sont satisfaites, on parlera de
\emph{quasi-géodésique de classe Morse} ou $(L,A,\epsilon,R)$-géodésique de classe Morse et on dira que le
groupe $\Gamma$ est \emph{quasi-plongé de classe Morse} ou
$(L,A,\epsilon,R)$-plongé de classe Morse dans
$\mathbf{SL}_n( \mathbf{R})$. Lorsqu'elles sont vérifiées seulement pour les
géodésiques de longueur $P\leq S$ donnée, on dira que $\Gamma$  est \emph{
  localement $(L,A,\epsilon,R, S)$-plongé de classe Morse} ou localement
quasi-plongé de classe Morse dans $\mathbf{SL}_n(
\mathbf{R})$. Puisque toutes les propriétés énoncées dans le théorème sont
invariantes par multiplication à gauche par un élément de $\Gamma$, on peut
aussi bien se restreindre aux géodésiques $( \gamma_p)_{p=0,\dots,P}$
vérifiant $\gamma_0 = e_\Gamma$. La condition d'être localement quasi-plongé
de classe Morse
n'implique donc qu'un nombre \emph{fini} de géodésiques (elles aussi finies)
de $\Gamma$. Cette condition suffit cependant à assurer la condition globale.

\begin{theo}[\cite{kapovichleebporti14}]
  \label{theo:local-global}
  Soit $(L,A,\epsilon,R)$ fixé et $\epsilon'> \epsilon$. Il existe $S\geq 0$
  tel que si $\Gamma$ est un sous-groupe de $\mathbf{SL}_n(\mathbf{R})$,
  hyperbolique au sens de Gromov, si $\Gamma$ est régulier, et si $\Gamma$ est
    localement $(L,A,\epsilon',R,S)$-plongé de classe Morse dans $\mathbf{SL}_n(
  \mathbf{R})$, alors $\Gamma$ est  $(L,A,\epsilon,R)$-plongé de classe Morse.
\end{theo}

\subsection{Divergence}
\label{sec:divergence}

Les représentations Anosov vérifient d'autres propriétés fortes de
divergence. Dans \textcite{ggkw_anosov} est prouvé le fait que, si $\Gamma$ est un
sous-groupe Anosov de $\mathbf{GL}_n( \mathbf{R})$, alors, pour tout
$i=2,\dots, n$ et pour tout rayon géodésique $(\gamma_p)_{p\in \mathbf{N}}$
dans $\Gamma$, l'application $\mathbf{N} \to \mathbf{R}_+ / p \mapsto \mu_1(
\gamma_p) - \mu_i(\gamma_p)$ est une quasi-isométrie. Les constantes de
quasi-isométries peuvent être prises uniformes pour les rayons géodésiques
tels que $\gamma_0=e_\Gamma$. Réciproquement il est établi que, si la fonction $\Gamma\to
\mathbf{R}_+ / \gamma \mapsto \mu_1(\gamma) - \mu_2(\gamma)$, pour un
sous-groupe $\Gamma$ de type fini, vérifie une
propriété de divergence \og logarithmique\fg{}, les suites $(f_p)_{p\in
  \mathbf{N}}$, dans $\mathcal{F}_{1,n-1}$, associées à un rayon géodésique
$(\gamma_p)_{p\in \mathbf{N}}$ par la relation $\gamma_p \cdot s_0 \in V( s_0,
f_p)$ sont toutes convergentes et de même limite ce qui permet de définir une
application $\beta: \partial_\infty \Gamma\to \mathcal{F}_{1,n-1}$
équivariante et continue. Lorsque de plus $\Gamma$ est hyperbolique au sens de
Gromov et que les fonctions $p \mapsto \mu_1(
\gamma_p) - \mu_i(\gamma_p)$ sont des quasi-isométries, il est démontré que $\beta$ vérifie la
propriété de transversalité et que le sous-groupe $\Gamma$ est Anosov.

Un résultat analogue est le fait que, pour toute quasi-géodésique de classe
Morse $(\gamma_p)_{p\in \mathbf{N}}$, la suite des projections de Cartan
 $( \boldsymbol{\mu}( \gamma_p \gamma_{0}^{-1}))_{p\in \mathbf{N}}$ est, elle
 aussi, une quasi-géodésique de classe Morse
 \parencite{kapovichleebporti17}.

\section{Hyperbolicité du groupe $\mathbf{\Gamma}$}
\label{sec:hyperb-du-groupe}

Cette section explique une partie des conclusions du théorème~\ref{th1}, à
savoir celle concernant l'hyperbolicité du groupe. Remarquons que les
hypothèses du théorème peuvent se réexprimer en disant que l'application
orbitale $\Gamma \to X_n / \gamma \mapsto \gamma\cdot s_0$ est un plongement
quasi-isométrique et qu'il existe $D\geq 0$ et $\epsilon>0$ tels que, pour
tous $\gamma$, $\gamma'$ dans $\Gamma$, si $d_{X_n}( \gamma\cdot s_0, \gamma'
\cdot s_0) \geq D$ alors le segment d'extrémités $\gamma\cdot s_0$ et $\gamma'
\cdot s_0$ est $\epsilon$-régulier. Posons $Z = \Gamma\cdot s_0 \subset X_n$
muni de la distance $d_Z$ déduite de celle de $X_n$ par restriction, $Z$ est
quasi-isométrique à $\Gamma$.

Il s'agit donc de démontrer que $(Z, d_Z)$ est hyperbolique au sens de
Gromov. La caractérisation suivante de l'hyperbolicité sera utilisée ici :
$(Z,d_Z)$ est hyperbolique au sens de Gromov si et seulement si tout cône
asymptotique (notion due à 
\textcite{dries_wilkie} et à
\textcite{gromov-asy-inv}) de $(Z,d_Z)$ est un arbre réel, i.e.\ est 0-hyperbolique.

Soit $\mathfrak{F} \subset \mathcal{P}( \mathbf{N})$ un ultrafiltre non
principal (i.e.\ $\mathfrak{F}$ ne contient aucun
ensemble fini, est stable
par intersection finie et par passage au sur-ensemble, et est maximal pour ces
propriétés). Un \emph{cône asymptotique}  de $(Z,d_Z)$ est un espace métrique $(
\mathcal{Z}, d_{\mathcal{Z}})$ obtenu à partir d'une suite $(\lambda_p)_{p\in
  \mathbf{N}}$ de réels strictement positifs tendant vers $0$ et d'une suite
de points base $(b_p)_{p\in \mathbf{N}} \in Z^{ \mathbf{N}}$ de la manière
suivante :
\begin{itemize}
\item soit $\widehat{ \mathcal{Z}}$ le sous-ensemble des suites $(z_p)_{p\in
    \mathbf{N}}$ telles que la suite réelle $( \lambda_p \, d_Z( b_p,
  z_p))_{p\in \mathbf{N}}$ est bornée selon $\mathfrak{F}$ (i.e.\ il existe
  $M\!\geq \!0$ tel que l'ensemble \mbox{$\{ p\!\in\! \mathbf{N} \!\mid\! | \lambda_p \, d_Z( b_p,
  z_p) |\! \leq \!M\}$} appartient à $\mathfrak{F}$). Dans ce cas, cette suite est
  automatiquement convergente selon $\mathfrak{F}$ (i.e.\ il existe $t\in
  \mathbf{R}$ tel que, pour tout $\epsilon>0$, $\{ p\!\in \!\mathbf{N} \mid |
  \lambda_p \, d_Z( b_p, z_p) -t | \leq \epsilon\}$ appartient à
  $\mathfrak{F}$)~;
\item on définit, quels que soient $(z_p)_{p\in \mathbf{N}}$ et $(z'_p)_{p\in
    \mathbf{N}}$ dans  $\widehat{ \mathcal{Z}}$, $d_{ \widehat{ \mathcal{Z}}}(
  (z_p)_p, (z'_p)_p) = \lim_{ \mathfrak{F}} \lambda_p\, d_Z( z_p,
  z'_p)$. Alors $d_{ \widehat{ \mathcal{Z}}}$ est symétrique et vérifie
  l'inégalité triangulaire~;
\item on obtient un espace métrique $\mathcal{Z}$ en quotientant $\widehat{
    \mathcal{Z}}$ par la relation d'équivalence dont le graphe est $\{ d_{
    \widehat{ 
      \mathcal{Z}}} =0\} \subset  \widehat{ \mathcal{Z}} \times  \widehat{
    \mathcal{Z}}$.
\end{itemize}

Ce cône asymptotique est maintenant un espace métrique connexe par arcs ; mieux il
existe des géodésiques (paramétrées par un intervalle de $\mathbf{R}$) entre
toute paire de points de $\mathcal{Z}$. 

Aussi (et presque par construction) $\mathcal{Z}$ se plonge dans le cône
asymptotique $\mathcal{X}$ de $X_n$. Par les travaux de 
\textcite{kleinerleeb-rig}, il est connu que ce cône asymptotique $\mathcal{X}$ a la
propriété remarquable d'être un immeuble euclidien de type $A_{n-1}$ (i.e.\ du
même type que $X_n$). Concrètement, dans $\mathcal{X}$, on peut définir :
\begin{itemize}
\item pour tout $(s,s')\in \mathcal{X}^2$, un $n$-uplet ${\boldsymbol{\mu}}(s,s') =
  (\mu_1, \dots, \mu_n) \in \mathbf{R}^n$ avec $\mu_1 \geq \mu_2 \geq \cdots
  \geq \mu_n$ (et $\sum_{i=1}^{n} \mu_i=0$) ; en particulier les notions de
  segment régulier et de segment \mbox{$\epsilon$-régulier}~;
  \item une \og variété drapeau\fg{} à l'infini $\mathcal{F}_\infty$ (ici l'ultraproduit $\prod_{
    \mathfrak{F}} \mathcal{F}_{1,n-1}$) et donc des cônes de Weyl et des
  diamants. Une différence étant que, pour les immeubles euclidiens, un
  segment régulier appartient à plusieurs cônes de Weyl, néanmoins les
  différents diamants que l'on peut dès lors construire sont tous égaux~;
\item en tout point $s$ de $\mathcal{X}$ une \og variété drapeau\fg{}
  $\mathcal{F}_s$ des directions des (germes de) cônes de Weyl au point $s$.
\end{itemize}
Dans ces variétés drapeaux $\mathcal{F}_\infty$ et $\mathcal{F}_s$, la notion
de transversalité est bien définie. Une différence importante avec le cas de
$X_n$ est que la topologie induite sur $\mathcal{F}_s$ est la topologie
discrète. Il y a des projections naturelles $\mathcal{F}_\infty \to
\mathcal{F}_s$ et, pour tout $s'$ dans $\mathcal{X}$ tel que le segment $ss'$
est régulier, le (germe de) cône de Weyl contenant ce segment définit un
élément $f(s') \in \mathcal{F}_s$. Cette fonction permet de caractériser le
diamant (de $\mathcal{X}$) $\diamondsuit_{s^-,s^+}$ : $s$ appartient au
diamant $\diamondsuit_{s^-,s^+}$ si et seulement si $f(s^-)$ et $f(s^+)\in
\mathcal{F}_s$ sont transverses. 

Par hypothèse sur $Z$ (ou plutôt sur $\Gamma$), l'ensemble $\mathcal{Z}$ est
automatiquement $\epsilon$-régulier, c'est-à-dire que toute paire de points
distincts dans $\mathcal{Z}\subset \mathcal{X}$ définit un segment
$\epsilon$-régulier. Une étape importante est d'obtenir :
\begin{quote}
  Tout segment rectifiable dans $\mathcal{Z}$ est contenu dans le diamant
  défini par ses extrémités, voir \cite[Th.~5.6]{kapovichleebporti14_2}.
\end{quote}
À partir de là, Kapovich, Leeb et Porti démontrent que $\mathcal{Z}$ est un
arbre métrique \parencite[Cor.~6.5 et 6.6]{kapovichleebporti14_2}.

Cette propriété des chemins rectifiables et $\epsilon$-réguliers est quant à
elle obtenue grâce à un contrôle des propriétés de contraction de la
projection $\pi_\diamondsuit \colon \mathcal{X} \to \diamondsuit$ sur un
diamant~$\diamondsuit$ (les diamants sont convexes). Cette projection est
$1$-lipschitzienne pour $d_{ \mathcal{X}}$ mais des estimées plus fines sont
nécessaires pour pouvoir conclure. Pour ce faire, Kapovich, Leeb et Porti
introduisent $d_\diamondsuit$ une métrique sur $\diamondsuit$ pour laquelle
ils démontrent :
\begin{quote}
  La projection $\pi_\diamondsuit$ est (localement) 1-lipschitzienne de
  $\mathcal{X} \setminus \diamondsuit$ dans $(  \diamondsuit, d_\diamondsuit
  )$ \parencite[Th.~4.8]{kapovichleebporti14_2}.
\end{quote}
La propriété recherchée provient alors du fait que, pour les segments
$\epsilon$-réguliers, la distance $d_\diamondsuit$ est quantitativement plus
grande que $d_\mathcal{X}$.

La définition de $d_\diamondsuit$ est ainsi : la longueur d'un chemin
géodésique par morceaux ne \og change\fg{} pas mais l'on autorise dans
$\diamondsuit = \diamondsuit_{s^-,s^+}$ uniquement des morceaux qui sont des
segments géodésiques $xy$ non inclus dans l'intérieur des cônes de Weyl $V(x,f)$ (ou
$V(y,f)$) qui contiennent $s^-$ ou $s^+$, autrement dit $y$ n'appartient pas
aux diamants $\diamondsuit_{s^-,x}$, $\diamondsuit_{x,s^+}$ \parencite[\S~4.1]{kapovichleebporti14_2}.

\section{Un panorama (trop rapide) autour des représentations Anosov}
\label{sec:un-panorama-trop}

Cet exposé n'a permis d'aborder que \emph{quelques} caractérisations des
représentations Anosov, les articles mentionnés plus haut 
en
contiennent d'autres renforçant encore le lien avec les sous-groupes
convexes-cocompacts (qui eux-mêmes admettent d'autres caractérisations).
En outre, il existe aussi des caractérisations utilisant les valeurs propres
des éléments du groupe plutôt que leurs valeurs principales.
Par
ailleurs, la métrique~$d_\diamondsuit$ introduite plus haut admet une
interprétation naturelle en termes de la géométrie finslérienne de l'espace
symétrique ; 
\textcite{kapovichleeb15} ont élaboré sur cette
géométrie finslérienne et les compactifications de l'espace symétrique qui
s'en suivent.

Il existe également une autre démonstration du lemme de Morse pour les espaces
symétriques donnée par 
\textcite{bochi_potrie_sambarino} et se basant sur la notion de décomposition
dominée en systèmes dynamiques et les résultats de 
\textcite{avila_bochi_yoccoz} et de 
\textcite{bochi_gourmelon} dans
ce domaine.

\subsection{Structures géométriques}
\label{sec:struct-geom}

Les représentations Anosov entretiennent des liens forts avec les structures
géométriques. Les résultats de Barbot sur les représentations Anosov dans
$\mathbf{PGL}_3( \mathbf{R})$ vont dans ce sens,
voir \S~\ref{sec:groupes-de-surfaces}. 
\textcite{guichard_wienhard_dod}
montrent, entre autres, que tout sous-groupe Anosov $\Gamma$ d'un groupe de Lie
semi-simple $G$ est l'holonomie d'une structure géométrique sur une variété
compacte ; plus précisément, il existe un $G$-espace homogène compact~$\mathcal{E}$, 
un ouvert $\Gamma$-invariant $\Omega$ de $\mathcal{E}$ (construit explicitement à
partir de l'application au bord $\partial_\infty \Gamma \to G/P$ et d'une
représentation de $G$ dans $\mathbf{GL}_N( \mathbf{R})$) sur lequel $\Gamma$
agit proprement et avec quotient compact. 
\textcite{kapovichleebportiv2}
offrent des éclairages nouveaux sur ces domaines de discontinuité avec quotient
compact : l'argument de cocompacité est de nature dynamique (et non cohomologique) et
$\Omega$ est ici un ouvert de $\mathcal{F}= G/P'$ d'une variété drapeaux et est
construit à partir d'une donnée combinatoire sur le groupe de Weyl $W$,
explicitement il s'agit d'un sous-ensemble de $W$ appelé \og épaississement
équilibré\fg{}. Y est démontré aussi que, dans $\mathcal{F}=
G/P_{\mathrm{min}}$, l'ouvert $\Omega$ est toujours non vide si $G$ a au moins
un facteur simple qui n'est pas de type $A_1$, $B_2$ ou $G_2$.

L'étude de ces structures géométriques prend maintenant un nouvel élan pour
les espaces de Teichmüller généralisés où les outils
analytico-algébro-géométriques, déjà utilisés depuis longtemps pour établir
des propriétés de ces espaces de modules, commencent aujourd'hui à être
employés pour étudier ces structures géométriques. Un premier papier dans
cette direction est celui de 
\textcite{collier_tholozan_toulisse} sur les représentations maximales dans
$\mathbf{SO}(2,n)$.

\subsection{Dynamique}
\label{sec:dynamique}

Le \og formalisme thermodynamique\fg{} a d'abord été introduit en théorie des
systèmes dynamiques par 
\textcite{bowen73,bowen}, 
\textcite{parry_pollicott}, 
\textcite{ruelle_book} et d'autres. Il a été
ensuite utilisé par 
\textcite{mcmullen_thermo} puis par 
\textcite{bridgeman_wp} pour les représentations fuchsiennes et quasi-fuchsiennes
(de nouvelles formules y sont données pour la métrique de Weyl et
Petersson). L'article de 
\textcite{bridgeman_canary_labourie_sambarino} développe ce formalisme pour les
représentations Anosov : ils construisent une nouvelle métrique (dite \og de
pression\fg{}) sur l'espace des modules des représentations Anosov (Zariski
denses) et démontrent que certains invariants dynamiques (entropie, dimension
de Hausdorff de l'ensemble limite, etc.) varient analytiquement. Avec ces
outils, 
\textcite{potrie_sambarino_pub} ont démontré un
résultat de rigidité pour les représentations $\rho: \Gamma_g \to
\mathbf{SL}_n( \mathbf{R})$ de la composante de Hitchin: la représentation~$\rho$ factorise par le
$\mathbf{SL}_2( \mathbf{R})$ principal si et seulement si l'entropie de $\rho$
est maximale (et donc égale à 1).

\subsection{Compactification, dégénérescence}
\label{sec:comp-degen}

Les travaux de Kapovich et Leeb et ceux de \textcite{gkw_tame} proposent des compactifications des espaces
localement symétriques $\Gamma \backslash G / K$ associés à des sous-groupes
Anosov. 
\textcite{kapovichleeb15} montrent que les
compactifications obtenues ont une structure naturelle de variétés à coins.
Réciproquement, ils établissent une caractérisation 
 à l'aide de ces compactifications : un sous-groupe $\Gamma$ de
$\mathbf{SL}_n( \mathbf{R})$ est Anosov si et seulement si $\Gamma$ est
uniformément régulier (c'est-à-dire les éléments de $\Gamma$ assez grands sont $\epsilon$-réguliers) et si le quotient $\Gamma\backslash X_n$ admet une compactification
finslérienne et \og respectant les fibrations naturelles à l'infini\fg{}.

Il est également possible
d'utiliser les structures géométriques mentionnées plus haut pour obtenir des
compactifications de quotients $\Gamma \backslash G / H$ plus
généraux. Centrons le reste de ce paragraphe sur les quotients $\Gamma
\backslash \left(G\times G\right)/ \Delta(G)$ où $\Gamma$ est un sous-groupe de $G\times G$
et $\Delta(G)$ est le sous-groupe diagonal de $G\times G$. L'étude de ces
quotients a une longue histoire que nous n'évoquons pas ici.
\textcite{ggkw_compact} montrent que ces doubles quotients ont des
compactifications naturelles lorsque~$\Gamma$ est Anosov dans un sur-groupe
$G'$ de $G\times G$ ($G'$ est explicite, il est simple si $G$ est simple).
Il y a des réciproques quand $G$ est
de rang réel égal à 1 : si $\Gamma$ est un sous-groupe de $G\times G$ agissant
proprement et avec quotient compact sur $(G\times G)/ \Delta(G)$, alors le
groupe $\Gamma$ est un sous-groupe Anosov dans $G'$ d'un groupe de Lie
simple $G'$ contenant $G\times G$ (voir \cite{ggkw_anosov} qui donne un énoncé
plus général).

Un autre aspect est la compactification de l'espace des modules des
représentations Anosov, et donc de comprendre les dégénérescences de suites de
représentations Anosov. Les travaux généraux de 
\textcite{parreau_publi}
construisent des compactifications de l'espace des modules de toutes les
représentations par des actions sur les immeubles euclidiens. 
\textcite{burgerpozzetti} 
analysent quels sont les immeubles
intervenant dans la compactification de l'espace des modules des
représentations maximales dans $\mathbf{Sp}_{2m}( \mathbf{R})$ et détaillent
les structures particulières des actions obtenues. 
\textcite{merlin}
donne quels groupes~$\Gamma$ admettent des suites \og fortement\fg{}
dégénérées d'actions Anosov.

L'action du groupe des automorphismes extérieurs $\mathrm{Out}( \Gamma)$ sur
l'espace des modules des représentations Anosov est propre
\parencite{labourie_energy,wienhard_mapping}. 
\textcite{CanaryLeeStover} 
introduisent une notion de représentation \og Anosov
amalgamés\fg{} et démontrent la propreté de l'action de $\mathrm{Out}( \Gamma)$
sur leur espace de modules quand le groupe $\Gamma$ a un bout. Ils montrent
également que, pour une infinité de groupes~$\Gamma$, l'espace des
représentations Anosov amalgamées, qui contient toujours l'espace des
représentations Anosov, ne coïncide pas avec ce dernier.

\subsection{Le retour des sous-groupes convexes-compacts}
\label{sec:le-retour-des}

Les travaux récents de 
\textcite{danciger_gueritaud_kassel_psdo-hyp,danciger_gueritaud_kassel-cc-proj}
rétablissent de plein droit la convexe-cocompacité dans l'étude des
représentations Anosov. L'un de leurs résultats est le théorème suivant : un
sous-groupe $\Gamma$ de $\mathbf{SL}_n( \mathbf{R})$ est Anosov et laisse
stable un convexe propre de l'espace projectif $\mathbb{P}^{n-1}( \mathbf{R})$
si et seulement si $\Gamma$ agit de manière convexe-cocompacte sur un ouvert
convexe $\Omega \subset \mathbb{P}^{n-1}( \mathbf{R})$ strictement convexe et
à bord $C^1$ (i.e.\ il existe $C\subset \Omega$ fermé, convexe,
$\Gamma$-invariant, sur lequel $\Gamma$ agit proprement avec quotient
compact). 
\textcite{zimmer-ccporj} démontre un résultat analogue pour une
notion légèrement différente de convexe cocompacité, et en déduit des énoncés
de rigidité. En fait des hypothèses plus faibles sur $\Omega$ impliquent déjà
que 
$\Gamma$ est Anosov. Si ce résultat ne peut s'appliquer aux sous-groupes de
$\mathbf{SL}_n(\mathbf{R})$ ne stabilisant aucun convexe de l'espace
projectif, on en tire tout de même le résultat inconditionnel suivant : un
sous-groupe $\Gamma$ de $\mathbf{SL}_n( \mathbf{R})$ est Anosov si et
seulement si l'action de $\Gamma$ par transconjugaison sur l'espace
$\mathbf{S}_n( \mathbf{R})$ des matrices symétriques est convexe-cocompacte :
il existe $\Omega \subset \mathbb{P} ( \mathbf{S}_n( \mathbf{R}))$, un ouvert convexe, $\Gamma$-invariant, strictement convexe et à bord
$C^1$, sur lequel
l'action de $\Gamma$ est convexe-cocompacte.

\bigskip


\printbibliography

\end{document}

